\documentclass{birkart}

 \newtheorem{thm}{Theorem}[section]
 \newtheorem{cor}[thm]{Corollary}
 \newtheorem{lem}[thm]{Lemma}
 \newtheorem{prop}[thm]{Proposition}
 \theoremstyle{definition}
 
 \theoremstyle{remark}
 \newtheorem{rem}[thm]{Remark}
 \newtheorem{prob}{Problem}
 \newtheorem{ex}{Example}
 \numberwithin{equation}{section}
 \usepackage{amsfonts}
\usepackage{amsmath}
\usepackage{amssymb}

\begin{document}
\title[On the {Bessmertny\u\i} Class]{On the {Bessmertny\u\i} Class of Homogeneous
Positive Holomorphic Functions of Several Variables}

 \author[Dmitry Kalyuzhny\u{\i}-Verbovetzki\u{\i}]{Dmitry S. Kalyuzhny\u{\i}-Verbovetzki\u{\i}}
\address{%
Department of Mathematics\\
The Weizmann Institute\\
Rehovot 76100\\
Israel}%
\email{dmitryk@wisdom.weizmann.ac.il}




\date{}

\begin{abstract}
The class of operator-valued functions which are homogeneous of
degree one, holomorphic in the open right polyhalfplane, have
positive semidefinite real parts there and take selfadjoint
operator values at real points, and its subclass consisting of
functions representable in the form of Schur complement of a block
of a linear pencil of operators with positive semidefinite
operator coefficients, are investigated. The latter subclass is a
generalization of the class of characteristic matrix functions of
passive 2n-poles considered as functions of impedances of its
elements, which was introduced by M.~F.~Bessmertny\u\i. Several
equivalent characterizations of the generalized Bessmertny\u{\i}
class are given, and its intimate connection with the Agler--Schur
class of holomorphic contractive operator-valued functions on the
unit polydisk is established.
\end{abstract}

\maketitle

\section{Introduction}\label{sec:intr}
In the Ph. D. Thesis of M.~F.~Bessmertny\u{\i} \cite{Bes}, which
appeared in Russian about twenty years ago and until very recent
time was unknown to Western readers (the translations of some its
parts into English are appearing now: see \cite{Bes1, Bes2,
Bes3}), rational $n\times n$ matrix-valued functions representable
in the form
\begin{equation}\label{eq:sc}
    f(z)=a(z)-b(z)d(z)^{-1}c(z),\quad z\in\mathbb{C}^N,
\end{equation}
with a linear $(n+p)\times(n+p)$ matrix-valued function
\begin{equation}\label{eq:lp}
    A(z)=A_0+z_1A_1+\cdots+z_NA_N=\left[\begin{array}{cc}
      a(z) & b(z) \\
      c(z) & d(z)
    \end{array}\right]
\end{equation}
were considered. Another form of such a representation is
\begin{equation}\label{eq:lr}
    f(z)=\left(\left[\begin{array}{cc}
      I_{n\times n} & 0_{n\times p} \end{array}\right]A(z)^{-1}
      \left[\begin{array}{c}
        I_{n\times n} \\
       0_{p\times n}\end{array}\right]\right)^{-1},
\end{equation}
and both (\ref{eq:sc}) and (\ref{eq:lr}) were called by
Bessmertny\u{\i} a \emph{long resolvent representation}. It is
easy to see that this is nothing but the \emph{Schur complement}
of the block $d(z)$ in the linear matrix pencil $A(z)$. In
\cite{Bes} (see also \cite{Bes1}) Bessmertny\u{\i} constructed a
long resolvent representation for an arbitrary rational matrix
function, and in the homogeneous case, i.e., when
\begin{equation}\label{eq:hom}
    f(\lambda z_1,\ldots,\lambda z_N)=\lambda
    f(z_1,\ldots,z_N),\quad \lambda\in\mathbb{C}\backslash\{ 0\},\
    z=(z_1,\ldots,z_N)\in\mathbb{C}^N,
\end{equation}
one has $A_0=0$.

A particular role in his thesis is played by functions of the form
(\ref{eq:sc}) or (\ref{eq:lr}) with $A_0=0$ and $A_k=A_k^T\geq 0,\
k=1,\ldots,N$ (matrices $A_k$ are assumed to have real entries),
due to their relation to electrical circuits. He proved that such
functions constitute the class (let us denote it by
$\mathbb{R}\mathcal{B}_N^{n\times n}$) of characteristic matrix
functions of passive 2n-poles, where impedances of elements
(resistances, capacitances, inductances and ideal transformers are
allowed) are considered as independent variables (let us note,
that in the analytic theory of electrical circuits it is customary
to consider characteristic matrices as functions of frequency,
e.g., see \cite{Kar, SeRe, EPo, LiFl}). It is easy to verify that
any $f\in\mathbb{R}\mathcal{B}_N^{n\times n}$ satisfies the
following properties:
\begin{equation}\label{eq:pos}
    f(z)+f(z)^*\geq 0,\quad z\in\Pi^N,
\end{equation}
where $\Pi^N:=\{ z\in\mathbb{C}^N:\ \mbox{Re}z_k>0,\
k=1,\ldots,N\}$ is the \emph{open right polyhalfplane}, i.e., the
Cartesian product of $N$ copies of the open right half-plane
$\Pi\subset\mathbb{C}$,
\begin{equation}\label{eq:real}
    f(\bar{z})=f(z)^*=f(\bar{z})^T,\quad z\in\mathbb{C}^N,
\end{equation}
where $\bar{z}:=(\overline{z_1},\ldots,\overline{z_N})$, together
with property (\ref{eq:hom}). Denote by
$\mathbb{R}\mathcal{P}_N^{n\times n}$ the class of rational
\emph{homogeneous positive real} $n\times n$ matrix functions,
i.e., rational functions  taking $n\times n$ matrix values and
satisfying conditions (\ref{eq:hom})--(\ref{eq:real}). Then
$\mathbb{R}\mathcal{B}_N^{n\times
n}\subset\mathbb{R}\mathcal{P}_N^{n\times n}$. Let us remark that
replacement of the requirement $A_k=A_k^T\geq 0$ by $A_k=A_K^*\geq
0,\ k=1,\ldots,N$ (i.e., removing the assumption that positive
semidefinite matrices $A_k$ have only real entries), in the
definition of $\mathbb{R}\mathcal{B}_N^{n\times n}$, and removing
the second equality in condition (\ref{eq:real}) in the definition
of $\mathbb{R}\mathcal{P}_N^{n\times n}$ define the classes
$\mathcal{B}_N^{n\times n}=\mathbb{C}\mathcal{B}_N^{n\times n}$
and $\mathcal{P}_N^{n\times n}=\mathbb{C}\mathcal{P}_N^{n\times
n}$, and $\mathcal{B}_N^{n\times n}\subset\mathcal{P}_N^{n\times
n}$. It is clear that in the case $N=1$ one has
$\mathbb{R}\mathcal{B}_N^{n\times
n}=\mathbb{R}\mathcal{P}_N^{n\times n}=\{ f(z)=zA:\ A=A^T\geq 0\}$
and $\mathcal{B}_N^{n\times n}=\mathcal{P}_N^{n\times n}=\{
f(z)=zA:\ A=A^*\geq 0\}$, where $z\in\mathbb{C}$ and $A$ is an
$n\times n$ matrix with real (resp., complex) entries, thus this
case is trivial. It was shown in \cite{Bes} that in the case $N=2$
one has $\mathbb{R}\mathcal{B}_N^{n\times
n}=\mathbb{R}\mathcal{P}_N^{n\times n}$, too (and we shall prove
in the present paper that $\mathcal{B}_N^{n\times
n}=\mathcal{P}_N^{n\times n}$). For $N\geq 3$ the question whether
the inclusion $\mathbb{R}\mathcal{B}_N^{n\times
n}\subset\mathbb{R}\mathcal{P}_N^{n\times n}$ (as well as
$\mathcal{B}_N^{n\times n}\subset\mathcal{P}_N^{n\times n}$) is
proper, is still open.

Bessmertny\u{\i} constructed the long resolvent representations
for the following special cases of functions from
$\mathbb{R}\mathcal{P}_N^{n\times n}$:
\begin{itemize}
    \item any $f\in\mathbb{R}\mathcal{P}_2^{n\times n}$;
    \item any rational scalar function $f=\frac{P}{Q}\in
    \mathbb{R}\mathcal{P}_N$ with co-prime polynomials $P$ and
    $Q$, where $\deg P=2$ (see \cite{Bes}, and also \cite{Bes3});
    \item any so-called \emph{primary} rational matrix function $f=\frac{P}{Q}\in
    \mathbb{R}\mathcal{P}_N^{n\times n}$, i.e., such that the
    matrix-valued polynomial $P$ and the scalar polynomial $Q$ are
    co-prime and of degree at most one with respect to each
    variable, and for each $j,k=1,\ldots,N$ there exist
    scalar polynomials
    $\varphi_{jk},\psi_1^{(j)}(z),\ldots,\psi_n^{(j)}(z)$ such
    that
    \begin{eqnarray*}
    \left(\frac{\partial Q}{\partial z_j}\right)\frac{\partial}{\partial z_k}
    \left(\frac{Q}{\frac{\partial Q}{\partial
    z_j}}\right)(z)=\varphi_{jk}^2(z);\\
    Q^2(z)\frac{\partial f}{\partial z_j}(z)=\left[\begin{array}{c}
      \psi_1^{(j)}(z) \\
      \vdots \\
      \psi_n^{(j)}(z)
    \end{array}\right]\left[\begin{array}{ccc}
      \psi_1^{(j)}(z) & \ldots & \psi_n^{(j)}(z)
    \end{array}\right];
    \end{eqnarray*}
it was shown that $f$ is primary if and only if there exists its
long resolvent representation where coefficients $A_k,\
k=1,\ldots,N$, of $A(z)$ in (\ref{eq:lp}) are positive
semidefinite (PSD) $(n+p)\times(n+p)$ matrices of rank one with
real entries, and $A_0=0$ (see \cite{Bes}, and also \cite{Bes2}).
\end{itemize}
However, no any inner characterization of the class
$\mathbb{R}\mathcal{B}_N^{n\times n}$ appears in \cite{Bes}, i.e.,
Bessmertnyi's thesis doesn't give an idea how to distinguish
functions from $\mathbb{R}\mathcal{P}_N^{n\times n}$ which admit a
long resolvent representation, except the cases mentioned above
and those functions which arise as characteristic functions of
certain concrete passive electrical 2n-poles.

The purpose of our paper is to give such a characterization. We
succeed, however in a more general framework, which is quite
natural. Firstly, we consider operator-valued functions instead of
only matrix-valued ones. Secondly, we consider holomorphic
functions instead of only rational ones. Thirdly, we permit
infinite-dimensional long resolvent representations instead of
only finite-dimensional ones.

We start with the ``complex case", i.e., generalize the classes
$\mathcal{B}_N^{n\times n}=\mathbb{C}\mathcal{B}_N^{n\times n}$
and $\mathcal{P}_N^{n\times n}=\mathbb{C}\mathcal{P}_N^{n\times
n}$. In Section~\ref{sec:classes} we introduce the corresponding
classes $\mathcal{B}_N(\mathcal{U})$ and
$\mathcal{P}_N(\mathcal{U})$ of homogeneous positive holomorphic
$L(\mathcal{U})$-valued functions (throughout this paper
$L(\mathcal{U, V})$ denotes the Banach space of bounded linear
operators mapping a Hilbert space $\mathcal{U}$ into a Hilbert
space $\mathcal{V}$, and $L(\mathcal{U}):=L\mathcal{(U, U)}$; all
Hilbert spaces are supposed to be complex). We obtain the
characterization of functions from the class
$\mathcal{B}_N(\mathcal{U})$ (which we call the
\emph{Bessmertny\u{\i} class}) via a couple of identities which
involve certain PSD kernels. In Section~\ref{sec:b-calc} we show
that one of these identities turns under the Cayley transform over
the variables into the Agler identity for holomorphic functions on
the \emph{unit polydisk} $\mathbb{D}^N:=\{ z\in\mathbb{C}^N:\
|z_k|<1,\ k=1,\ldots,N\}$ taking operator values with PSD real
parts. The latter means that the image of the
$\mathcal{B}_N(\mathcal{U})$ under the Cayley transform over the
variables is a subclass in the Agler--Herglotz class
$\mathcal{AH}_N(\mathcal{U})$, introduced in \cite{Ag}. Using the
characterization of $\mathcal{AH}_N(\mathcal{U})$ in terms of
functional calculus of $N$-tuples of commuting strictly
contractive linear operators on a Hilbert space, we obtain the
characterization of $\mathcal{B}_N(\mathcal{U})$ in terms of
functional calculus of $N$-tuples of commuting bounded strictly
accretive operators on a Hilbert space. In Section~\ref{sec:Ag} we
characterize the image of the Bessmertny\u{\i} class
$\mathcal{B}_N(\mathcal{U})$ under the double Cayley transform
(``double" means that this linear-fractional transform is applied
to the variables and to the operator values, simultaneously), as a
subclass in the Agler--Schur class $\mathcal{AS}_N(\mathcal{U})$,
also introduced in \cite{Ag}. This characterization turns out to
be pretty surprising: a function belongs to this subclass if and
only if it is representable as a transfer function of an Agler
unitary colligation for which the colligation operator is not only
unitary, but also selfadjoint. In Section~\ref{sec:dehom} we
establish a natural one-to-one correspondence between
$\mathcal{B}_N(\mathcal{U})$ (which consists of homogeneous
functions of $N$ variables) and certain class of (generically)
non-homogeneous functions of $N-1$ variables. However, a special
complicated structure of the latter class does rather convince us
that $\mathcal{B}_N(\mathcal{U})$ is more likeable to deal with.
In Section~\ref{sec:real} we turn to the ``real" case. We
introduce the notions of $\iota$-real operator and $\iota$-real
operator-valued function for an anti-unitary involution $\iota
=\iota_{\mathcal{U}}$ on a Hilbert space $\mathcal{U}$ which plays
a role analogous to the complex conjugation in $\mathbb{C}$, and
then introduce the subclasses
$\iota\mathbb{R}\mathcal{B}_N(\mathcal{U})$ and
$\iota\mathbb{R}\mathcal{P}_N(\mathcal{U})$ in the classes
$\mathcal{B}_N(\mathcal{U})$ and $\mathcal{P}_N(\mathcal{U})$,
respectively, consisting of $\iota$-real operator-valued
functions. These subclasses generalize the classes
$\mathbb{R}\mathcal{B}_N^{n\times n}$ and
$\mathbb{R}\mathcal{P}_N^{n\times n}$, respectively. We adapt the
results for the ``complex case" of the preceding sections to this
``real case", i.e., give the characterizations of
$\iota\mathbb{R}\mathcal{B}_N(\mathcal{U})$ in terms of long
resolvent representations, in terms of identities involving PSD
kernels, in terms of Agler's unitary colligations and their
transfer function representations for images of its elements under
the double Cayley transform. In Section~\ref{sec:conclusion} we
summarize the results obtained in this paper, and also formulate
and briefly discuss the most important open problems arising in
connection with our investigation.

\section{The classes of homogeneous positive holomorphic
functions}\label{sec:classes}
\subsection{}
Let $\mathcal{U}$ be a Hilbert space. Consider the class
$\mathcal{P}_N(\mathcal{U})$ consisting of all
$L\mathcal{(U)}$-valued functions $f$ holomorphic in the domain
$\Omega_N:=\bigcup_{\lambda\in\mathbb{T}}(\lambda\Pi)^N\subset\mathbb{C}^N$
(here for a fixed $\lambda\in\mathbb{T}$ a polyhalfplane
$(\lambda\Pi)^N$ is the product of $N$ copies of the half-plane
$\lambda\Pi:=\{\lambda z=(\lambda z_1,\ldots,\lambda z_N):\
z\in\Pi^N\}$), such that the following conditions are satisfied:
\begin{eqnarray}
f(\lambda z_1,\ldots,z_N)=\lambda f(z_1,\ldots,z_N), &
\lambda\in\mathbb{C}\backslash\{ 0\},\ z\in\Omega_N;
\label{eq:ghom}\\
f(z)+f(z)^*\geq 0, & z\in\Pi^N; \label{eq:gpos}\\
f(\bar{z})=f(z)^*, & z\in\Omega_N. \label{eq:gsym}
\end{eqnarray}

Let us formulate and prove some geometrical properties of the
domain $\Omega_N$ which appears naturally in this definition of
the class $\mathcal{P}_N(\mathcal{U})$, and the definition of the
class $\mathcal{B}_N(\mathcal{U})$ given subsequently (see
Remark~\ref{rem:dom} below for the motivation), even though we
will not use these properties explicitly in this paper.
\begin{prop}
The domain $\Omega_N$ has the following properties:
\begin{description}
    \item[(i)] $\Omega_N$ is an (open) \textbf{cone} in $\mathbb{C}^N$, i.e., for
each $z=(z_1,\ldots,z_N)\in\Omega_N$ and any real $t>0$ one has
$tz=(tz_1,\ldots,tz_N)\in\Omega_N$;
    \item[(ii)] $\Omega_N$ is a \textbf{circular domain}, i.e., for each
$z=(z_1,\ldots,z_N)\in\Omega_N$ and any $\lambda\in\mathbb{T}^N$
one has $\lambda z=(\lambda z_1,\ldots,\lambda z_N)\in\Omega_N$;
    \item[(iii)] for $N=1$ and $N=2$ one has
$\mbox{clos}(\Omega_N)=\mathbb{C}^N$, and for $N\geq 3$ one has
$\mbox{clos}(\Omega_N)\neq\mathbb{C}^N$;
    \item[(iv)] $\mbox{clos}(\Omega_N)$, as well as $\Omega_N$, is
    not convex;
    \item[(v)] $\Omega_N$ is \textbf{pseudo-convex} (i.e., there exists
    a plurisubharmonic function on $\Omega_N$, which tends to $+\infty$ as its
    variable point approaches to the boundary $\partial\Omega_N$);
    or equivalently,
    \item[(vi)] $\Omega_N$ is a \textbf{holomorphy domain} (i.e., there exists a
    holomorphic function on $\Omega_N$, which is not holomorphically extendable
    to any bigger domain); or equivalently,
    \item[(vii)] $\Omega_N$ is \textbf{not holomorphically
    extendable to any boundary point} (i.e., for any point $a\in\partial\Omega_N$
    there exist a neighborhood $\Gamma$ of $a$, and a function
    which is holomorphic in $\Omega_N\cap\Gamma$ and not holomorphically
    extendable to $a$).
\end{description}
For the proof of equivalence of properties \textbf{(v)--(vii)},
see \cite{Sh2}.
\end{prop}
\begin{proof}
The properties (i) and (ii) are evident.

(iii). $N=1$:
$\mbox{clos}(\Omega_1)=\mbox{clos}(\mathbb{C}\backslash\{
0\})=\mathbb{C}$.

$N=2$: for any $z=(z_1,z_2)\in\mathbb{C}^2$ there is a
$\lambda\in\mathbb{T}$ such that $z_1\in\mbox{clos}(\lambda\Pi),
z_2\in\mbox{clos}(\lambda\Pi)$, thus $z=(z_1,z_2)\in
(\mbox{clos}(\lambda\Pi))^2\subset\mbox{clos}(\Omega_2)$, i.e.,
$\mbox{clos}(\Omega_2)=\mathbb{C}^2$.

$N\geq 3$: a point $z=(1,\exp(2\pi i/N),\ldots,\exp(2\pi
i(N-1)/N))$ doesn't belong, together with some neighborhood, to
any open polyhalfplane, and therefore to $\Omega_N$, thus
$z\notin\mbox{clos}(\Omega_N)$.

(iv). Clearly, $\Omega_N$ is not convex because for any
$z\in\Omega_N$ one has $-z\in\Omega_N$, and
$\frac{z+(-z)}{2}=0\notin\Omega_N$.

Let us show that
$\mbox{clos}(\Omega_N)=\bigcup_{\lambda\in\mathbb{T}}(\lambda\mbox{clos}(\Pi))^N$.
Indeed, if $z\in(\lambda\mbox{clos}(\Pi))^N$ for some
$\lambda\in\mathbb{T}$ then, clearly, $z\in\mbox{clos}(\Omega_N)$.
Conversely, if
$z\notin\bigcup_{\lambda\in\mathbb{T}}(\lambda\mbox{clos}(\Pi))^N$
then for any $\lambda\in\mathbb{T}$ one has
$\frac{z}{\lambda}=(\frac{z_1}{\lambda},\ldots,\frac{z_N}{\lambda})\notin
(\mbox{clos}(\Pi))^N$, moreover, there exist a neighborhood
$\Gamma_\lambda\subset\mathbb{C}^N$ of $z$, and a neighborhood
$\Delta_\lambda\subset\mathbb{C}$ of $\lambda$ such that
 for any $z'\in\Gamma_\lambda$ and
$\lambda'\in\Delta_\lambda$ one has $\frac{z'}{\lambda'}\notin
(\mbox{clos}(\Pi))^N$. Since the collection of open sets
$\{\Delta_\lambda\}_{\lambda\in\mathbb{T}}$ is a covering of a
compact set $\mathbb{T}$, one may choose its finite subcovering
$\{\Delta_{\lambda_j}\}_{j=1}^m$. Set
$\Gamma:=\bigcap_{j=1}^m\Gamma_{\lambda_j}$. Then for any
$z'\in\Gamma$ and $\lambda\in\mathbb{T}$ one has
$\frac{z'}{\lambda}\notin (\mbox{clos}(\Pi))^N$, i.e., $\Gamma\cap
(\lambda\mbox{clos}(\Pi))^N=\emptyset$. Thus,
$\Gamma\cap\Omega_N=\emptyset$, i.e.,
$z\notin\mbox{clos}(\Omega_N)$.

A point $z=(1,\exp(2\pi i/N),\ldots,\exp(2\pi
i(N-1)/N))\notin\mbox{clos}(\Omega_N)$ can be expressed as
$z=\frac{1}{N}\sum_{k=1}^Nz^{(k)}$, where
$$z^{(k)}=(0,\ldots,0,\underbrace{N\exp(2\pi i(k-1)/N)}_{k-\mbox{th place}},0,\ldots,0)
\in\bigcup_{\lambda\in\mathbb{T}}(\lambda\mbox{clos}(\Pi))^N=\mbox{clos}(\Omega_N),$$
thus the set $\mbox{clos}(\Omega_N)$ is not convex.

(vii). First of all, let us show that the boundary points of
$\Omega_N$ can be of two different types:
\begin{enumerate}
    \item
    $z^\circ=(z_1^\circ,\ldots,z_{j-1}^\circ,0,z_{j+1}^\circ,
    \ldots,z_N^\circ)$ for some $j\in\{ 1,\ldots,N\}$, where ${z_j^\circ}':=
    (z_1^\circ,\ldots,z_{j-1}^\circ,z_{j+1}^\circ,
    \ldots,z_N^\circ)\in\mbox{clos}(\Omega_{N-1})$;
    \item $z^\circ=(z_1^\circ,\ldots,z_N^\circ)\in\mathbb{C}^N$ such that $z_k^\circ\neq
    0$ for all $k\in\{ 1,\ldots,N\}$, and there exist $l,j\in\{ 1,\ldots,N\},\ l<j$,
     and $t>0$ for which $z_l^\circ+tz_j^\circ=0$ and ${z_{lj}^\circ }'':=
    (z_1^\circ,\ldots,z_{l-1}^\circ,z_{l+1}^\circ,\ldots,
    z_{j-1}^\circ,z_{j+1}^\circ,
    \ldots,z_N^\circ)\in (e^{i\arg
    z_l^\circ}\mbox{clos}(\Pi))^{N-2}$ (or $\in (e^{-i\arg
    z_l^\circ}\mbox{clos}(\Pi))^{N-2}$).
\end{enumerate}
It is clear that points of these two types are boundary for
$\Omega_N$. Assume $z^\circ\in\partial\Omega_N$ is neither of
these types. Then $z^\circ\in(\lambda\mbox{clos}(\Pi))^N$ for some
$\lambda\in\mathbb{T}$, moreover
$z_{j_\mu}\in\partial(\lambda\Pi),\ \mu=1,\ldots,s$, for some
$j_1,\ldots,j_s\in\{ 1,\ldots,N\}$, $$\arg z_{j_1}=\ldots=\arg
z_{j_s}=\arg\lambda+\frac{\pi}{2}\ \left(\mbox{or}\
=\arg\lambda-\frac{\pi}{2}\right),$$ and $z_k\in(\lambda\Pi)$ for
$k\in\{ 1,\ldots,N\}\backslash\{ j_1,\ldots,j_s\}$. In this case,
there exists $\varepsilon>0$ such that $z^\circ\in(\lambda
e^{i\varepsilon}\Pi)^N$ (resp., $z^\circ\in(\lambda
e^{-i\varepsilon}\Pi)^N$), i.e., $z^\circ$ is an inner point of
$\Omega_N$, that contradicts to our assumption.

Let $\Gamma$ be a neighborhood of a point $z^\circ$ of type 1. Set
$f(z):=\frac{1}{z_j}$. Since
    for all points of $\Omega_N$ one has $z_j\neq 0$, $f$ is
    holomorphic in $\Omega_N\cap\Gamma$, and not holomorphically
    extendable to $z^\circ$.

Let $\Gamma$ be a neighborhood of a point $z^\circ$ of type 2. Set
$f(z):=\frac{1}{z_l+tz_j}$. Since
    for all points of $\Omega_N$ one has $z_l+tz_j\neq 0$, $f$ is
    holomorphic in $\Omega_N\cap\Gamma$, and not holomorphically
    extendable to $z^\circ$.

The proof is complete.
\end{proof}
An equivalent definition of the class $\mathcal{P}_N(\mathcal{U})$
will be given below. For that let us prove, first of all, the
following lemma.
\begin{lem}\label{lem:eqpos}
Let a scalar-valued function $f$ be holomorphic on
$\mathbb{C}\backslash\{ 0\}$ and satisfy $f(i^n\Pi)\subset i^n\Pi$
for $n=0,1,2,3$. Then f(z)=$\alpha z,\ z\in\mathbb{C}$, with some
number $\alpha >0$.
\end{lem}
\begin{proof}
Set $h(z):=\frac{f(z)}{z}$. From the assumptions of this lemma it
follows that $h$ doesn't take values in the negative real
semi-axis. Hence, the function $g(z):=\sqrt{h(z)}$ (with the
principal branch of the square root) is well defined and
holomorphic on $\mathbb{C}\backslash\{ 0\}$, and
$g(\mathbb{C}\backslash\{ 0\})\subset\Pi$. Since for any number
$w$ in the open left half-plane there is no sequence ${z_j}$ such
that $\lim_{j\rightarrow\infty}z_j=0$ and
$\lim_{j\rightarrow\infty}g(z_j)=w$, the point $z=0$ can not be
essentially singular for $g$ (by the Sokhotsky theorem, see e.g.
\cite{Sh1}). Hence, this point is not essentially singular for
$h$, too. But $z=0$ can not be a pole of $h$, since in this case
$\arg h(z)$ would take all values from $[-\pi, \pi)$, including
the value $-\pi$, which is banned. Thus, $h$ has a removable
singularity at $z=0$.

Now consider the function
$\widetilde{h}(z):=h(\frac{1}{z})=zf(\frac{1}{z})$. This function
is holomorphic on $\mathbb{C}\backslash\{ 0\}$, and also doesn't
take values in the negative real semi-axis. Applying the same
reasoning, we obtain that $\widetilde{h}$ has a removable
singularity at $z=0$, and therefore $h$ has a removable
singularity at infinity. Finally, we have got the entire function
$h$ which has no singularity at infinity. Thus, by the Liouville
theorem, $h$ is a constant. Therefore, $f(z)=\alpha z$, and since
$f(\Pi)=\alpha\Pi\subset\Pi$, we get $\alpha >0$.
\end{proof}
\begin{cor}
For an $L\mathcal{(U)}$-valued function $f$ holomorphic on the
domain $\Omega_N$, conditions (\ref{eq:ghom}) and (\ref{eq:gpos})
together are equivalent to the following set of conditions:
\begin{equation}\label{eq:eqpos}
    \begin{array}{cccc}
      f(z)+f(z)^*\geq 0, & z\in\Pi^N, \\
      f(z)+f(z)^*\leq 0, & z\in(-\Pi)^N, \\
      i(f(z)^*-f(z))\geq 0, & z\in(i\Pi)^N, \\
      i(f(z)^*-f(z))\leq 0, & z\in(-i\Pi)^N.
    \end{array}
\end{equation}
\end{cor}
\begin{proof}
If conditions (\ref{eq:ghom}) and (\ref{eq:gpos}) are fulfilled,
then (\ref{eq:eqpos}) follows immediately. Conversely, let
conditions (\ref{eq:eqpos}) are fulfilled. Consider for each fixed
$x\in\Pi^N\cap\mathbb{R}^N$ the cut-function
$\varphi_x(\lambda):=f(\lambda x_1,\ldots,\lambda x_N),\
\lambda\in\mathbb{C}\backslash\{ 0\}$, and for each fixed
$x\in\Pi^N\cap\mathbb{R}^N,\ u\in\mathcal{U}$ the function
$\psi_{x,u}(\lambda)=\langle\varphi_x(\lambda)u,u\rangle,\
\lambda\in\mathbb{C}\backslash\{ 0\}$. If
$\psi_{x,u}(\lambda)\equiv 0$ set $\alpha_{x,u}:=0$. If
$\psi_{x,u}(\lambda)\not\equiv 0$ then $\psi_{x,u}$ satisfies the
conditions of Lemma~\ref{lem:eqpos}. Thus,
$\psi_{x,u}(\lambda)=\alpha_{x,u}\lambda,\
\lambda\in\mathbb{C}\backslash\{ 0\}$ for some number
$\alpha_{x,u}>0$. Since $\psi_{x,u}(1)=\alpha_{x,u}=\langle
f(x)u,u\rangle$, one has
$\psi_{x,u}(\lambda)=\langle\varphi_x(\lambda)u,u\rangle=\langle
f(x)u,u\rangle\lambda$ for any fixed $x\in\Pi^N\cap\mathbb{R}^N,\
u\in\mathcal{U}$. Therefore, $\varphi_x(\lambda)=\lambda f(x)$,
i.e., $f(\lambda x_1,\ldots,\lambda x_N)=\lambda
f(x_1,\ldots,x_N),\ \lambda\in\mathbb{C}\backslash\{ 0\}$ for any
$x\in\Pi^N\cap\mathbb{R}^N$. By the uniqueness theorem for
holomorphic functions of several variables (see, e.g.,
\cite{Sh2}), we obtain (\ref{eq:ghom}). Since (\ref{eq:gpos}) is
contained in the set of conditions (\ref{eq:eqpos}), the proof is
complete.
\end{proof}
Finally, we get the following theorem.
\begin{thm}
An $L\mathcal{(U)}$-valued holomorphic function $f$ on $\Omega_N$
belongs to the class $\mathcal{P}_N(\mathcal{U})$ if and only if
conditions (\ref{eq:eqpos}) and (\ref{eq:gsym}) are satisfied.
\end{thm}

\subsection{}
Let $\mathcal{U}$ be a Hilbert space. Define the
\emph{Bessmertny\u{\i} class} $\mathcal{B}_N(\mathcal{U})$ as a
class of $L\mathcal{(U)}$-valued functions $f$ holomorphic on the
domain $\Omega_N$ and representable in the form
\begin{equation}\label{eq:gsc}
    f(z)=a(z)-b(z)d(z)^{-1}c(z),\quad z\in\Omega_N,
\end{equation}
where
\begin{equation}\label{eq:glp}
    A(z)=z_1A_1+\cdots +z_NA_N=\left[\begin{array}{cc}
      a(z) & b(z) \\
      c(z) & d(z)
    \end{array}\right]\in L(\mathcal{U\oplus H}),\quad z\in\Omega_N,
\end{equation}
for some Hilbert space $\mathcal{H}$, and bounded linear operators
$A_k=A_k^*,\ k=1,\ldots,N,$ on $\mathcal{U\oplus H}$ are PSD. It
is clear that $c(z)=\widetilde{b}(z):=b(\overline{z})^*$, i.e.,
$z_1c_1+\cdots +z_Nc_N=z_1b_1^*+\cdots +z_Nb_N^*,\ z\in\Omega_N$,
and linear pencils $a(z)=z_1a_1+\cdots +z_Na_N=\widetilde{a}(z),\
d(z)=z_1d_1+\cdots +z_Nd_N=\widetilde{d}(z)$ have PSD coefficients
$a_k=a_k^*,\ d_k=d_k^*,\ k=1,\ldots,N$, from $L\mathcal{(U)}$ and
$L\mathcal{(H)}$, respectively. For a function
$f\in\mathcal{B}_N(\mathcal{U})$ it is easy to check properties
(\ref{eq:ghom})--(\ref{eq:gsym}), thus
$\mathcal{B}_N(\mathcal{U})\subset\mathcal{P}_N(\mathcal{U})$.
\begin{rem}\label{rem:dom}
If a function $f$ is holomorphic on $\Pi^N$ and has the
representation (\ref{eq:gsc}) there, then $f$ can be extended to
$\Omega_N$ by homogeneity of degree one, and this extension is,
clearly, holomorphic and has a representation (\ref{eq:gsc}) in
$\Omega_N$. That is why we define the class
$\mathcal{B}_N(\mathcal{U})$ straight away as a class of functions
on $\Omega_N$. Keeping in mind the possibility and uniqueness of
such extension, we will write sometimes
$f\in\mathcal{B}_N(\mathcal{U})$ for functions defined originally
on $\Pi^N$.
\end{rem}
\begin{rem}\label{rem:glr}
If $f\in\mathcal{B}_N(\mathcal{U})$ and $f(z)$ is boundedly
invertible in $\Pi^N$ (and hence in $\Omega_N$) then
(\ref{eq:gsc}) can be rewritten in the form
\begin{equation}\label{eq:glr}
    f(z)=\left(P_\mathcal{U}A(z)^{-1}|\mathcal{U}\right)^{-1},\quad
    z\in\Omega_N,
\end{equation}
where $P_\mathcal{U}$ is the orthogonal projector onto
$\mathcal{U}$ in $\mathcal{U\oplus H}$, and $A(z)$ is given by
(\ref{eq:glp}). This follows from the equality
\begin{eqnarray}\nonumber\lefteqn{\left[\begin{array}{cc}
 a(z)-b(z)d(z)^{-1}c(z) & 0 \\
  0 & d(z)\end{array}\right]=}  \\ & & \left[\begin{array}{cc}
    I_\mathcal{U} & -b(z)d(z)^{-1} \\
    0 & I_\mathcal{H}
  \end{array}\right]\left[\begin{array}{cc}
      a(z) & b(z) \\
      c(z) & d(z)
    \end{array}\right]\left[\begin{array}{cc}
    I_\mathcal{U} & 0 \\
    -d(z)^{-1}c(z) & I_\mathcal{H}
  \end{array}\right]. \label{eq:schur-com} \end{eqnarray}
\end{rem}
Let $\Omega$ be a set. Recall that an $L\mathcal{(U)}$-valued
function $\Phi$ on $\Omega\times\Omega$ is called a
\emph{Hermitian symmetric positive semidefinite kernel} (or, for
the shortness, a \emph{PSD kernel}) if
\begin{equation}\label{eq:K1}
 \Phi(z,\zeta)=\Phi(\zeta,z)^*,\quad
(z,\zeta)\in\Omega\times\Omega,
\end{equation}
and for any $m\in\mathbb{N},\ \{
z^{(\mu)}\}_{\mu=1}^m\subset\Omega,\ \{
u^{(\mu)}\}_{\mu=1}^m\subset\mathcal{U}$ one has
\begin{equation}\label{eq:K2}
\sum\limits_{\mu=1}^m\sum\limits_{\nu=1}^m
\langle\Phi(z^{(\mu)},z^{(\nu)})u^{(\mu)},u^{(\nu)}\rangle\geq 0.
\end{equation}
If $\Omega$ is a domain in $\mathbb{C}^N$ and a PSD kernel
$\Phi(z,\zeta)$ on $\Omega\times\Omega$ is holomorphic in $z$ and
anti-holomorphic in $\zeta$, then $\Phi(z,\zeta)$ is said to be a
\emph{holomorphic PSD kernel on $\Omega\times\Omega$} (not to be
confused with a holomorphic function on $\Omega\times\Omega$).
\begin{thm}\label{thm:ker}
Let $f$ be an $L\mathcal{(U)}$-valued function holomorphic on
$\Pi^N$. Then $f\in\mathcal{B}_N(\mathcal{U})$ if and only if
there exist holomorphic PSD kernels $\Phi_k(z,\zeta),\
k=1,\ldots,N$, on $\Pi^N\times\Pi^N$ such that
\begin{equation}\label{eq:kerid}
    f(z)=\sum\limits_{k=1}^Nz_k\Phi_k(z,\zeta),\quad
    (z,\zeta)\in\Pi^N\times\Pi^N,
\end{equation}
holds. In this case the kernels $\Phi_k(z,\zeta),\ k=1,\ldots,N$,
can be uniquely extended to the holomorphic PSD kernels on
$\Omega_N\times\Omega_N$ (we denote the extended kernels by the
same letters) which are homogeneous of degree zero, i.e., for any
$\lambda\in\mathbb{C}\backslash\{ 0\}$
\begin{equation}\label{eq:homzero}
    \Phi_k(\lambda z,\lambda\zeta)=\Phi_k(z,\zeta),\quad
    (z,\zeta)\in\Omega_N\times\Omega_N,
\end{equation}
 and identity (\ref{eq:kerid}) is extended to all
$(z,\zeta)\in\Omega_N\times\Omega_N$.
\end{thm}
\begin{proof}
\textbf{Necessity.} This part of the theorem was proved by
Bessmertny\u{\i} in \cite{Bes} (see also \cite{Bes1}) for rational
matrix-valued functions $f$ by playing with a long resolvent
representation of $f$ in the form (\ref{eq:lr}). We follow the
same idea, however we use in our proof the representation
(\ref{eq:gsc}), which does exist for our case without an
additional assumption on the invertibility of $f(z)$ (see
Remark~\ref{rem:glr}).

Let $f\in\mathcal{B}_N(\mathcal{U})$. Then (\ref{eq:gsc}) holds
for some Hilbert space $\mathcal{H}$ and linear pencil of
operators $A(z)$ of the form (\ref{eq:glp}). Define
$$\psi(z):=\left[\begin{array}{c}
  I_\mathcal{U} \\
  -d(z)^{-1}c(z)
\end{array}\right]\in L\mathcal{(U,U\oplus H)},\quad z\in\Omega_N.$$
Then for all $(z,\zeta)\in\Omega_N\times\Omega_N$ one has
\begin{eqnarray*}
f(z) &=& a(z)-b(z)d(z)^{-1}c(z)\\
     &=& \left[\begin{array}{cc}
       I_\mathcal{U} & -c(\zeta)^*d(\zeta)^{*-1}
     \end{array}\right]\left[\begin{array}{c}
       a(z)-b(z)d(z)^{-1}c(z) \\
       0
     \end{array}\right]\\
     &=& \left[\begin{array}{cc}
       I_\mathcal{U} & -c(\zeta)^*d(\zeta)^{*-1}
     \end{array}\right]\left[\begin{array}{cc}
       a(z) & b(z) \\
       c(z) & d(z)
     \end{array}\right]\left[\begin{array}{c}
       I_\mathcal{U} \\
       -d(z)^{-1}c(z)
     \end{array}\right]\\
     &=& \psi(\zeta)^*A(z)\psi(z).
\end{eqnarray*}
Set $\Phi_k(z,\zeta):=\psi(\zeta)^*A_k\psi(z),\ k=1,\ldots,N$.
Clearly, functions $\Phi_k(z,\zeta),\ k=1,\ldots,N,$ satisfy
properties (\ref{eq:K1}) and (\ref{eq:K2}) for $\Omega=\Omega_N$.
Since $\psi$ is holomorphic on $\Omega_N$ we obtain that
$\Phi_k(z,\zeta),\ k=1,\ldots,N$, are holomorphic PSD kernels on
$\Omega_N\times\Omega_N$. Rewriting the equality
$$f(z)=\psi(\zeta)^*A(z)\psi(z),\quad (z,\zeta)\in\Omega_N\times\Omega_N,$$
in the form
\begin{equation}\label{eq:b}
    f(z)=\sum_{k=1}^Nz_k\Phi_k(z,\zeta),\quad (z,\zeta)\in\Omega_N\times\Omega_N,
\end{equation}
we obtain, in particular, (\ref{eq:kerid}).

\textbf{Sufficiency.} Let us note that the assumption that a
Hilbert space $\mathcal{H}$ involved implicitly in the
representation (\ref{eq:gsc}) can be infinite-dimensional, is
essential in our proof of this part of the theorem
(Bessmertny\u{\i} proved only the necessity part, under his
assumptions).

Let $f$ be an $L\mathcal{(U)}$-valued function holomorphic on
$\Pi^N$ and representable there in the form (\ref{eq:kerid}) with
some holomorphic PSD kernels $\Phi_k(z,\zeta),\ k=1,\ldots,N$. For
these kernels there exist auxiliary Hilbert spaces $\mathcal{M}_k$
and holomorphic $L\mathcal{(U,M}_k)$-valued functions $\varphi_k$
on $\Pi^N$ such that
$\Phi_k(z,\zeta)=\varphi_k(\zeta)^*\varphi_k(z),\
(z,\zeta)\in\Pi^N\times\Pi^N$ (see \cite{Aro}). Set
$\mathcal{M}:=\bigoplus_{k=1}^N\mathcal{M}_k,\
P_k:=P_{\mathcal{M}_k},\
\varphi(z):=\mbox{col}\left[\begin{array}{ccc}
  \varphi_1(z) & \ldots & \varphi_N(z)
\end{array}\right]\in L\mathcal{(U,M)},\ e:=(1,\ldots,1)\in\Pi^N$.
>From (\ref{eq:kerid}) we get
\begin{equation}\label{eq:zeta-e}
    f(e)=\sum_{k=1}^N\varphi_k(\zeta)^*\varphi_k(e),\quad
    \zeta\in\Pi^N.
\end{equation}
In particular,
\begin{equation}\label{eq:e-e}
    f(e)=\sum_{k=1}^N\varphi_k(e)^*\varphi_k(e).
\end{equation}
By subtracting (\ref{eq:e-e}) from (\ref{eq:zeta-e}) we get
$$\sum_{k=1}^N[\varphi_k(\zeta)-\varphi_k(e)]^*\varphi_k(e)=0,\quad
    \zeta\in\Pi^N,$$
i.e., the following orthogonality relation holds:
$$\mathcal{H}:=\mbox{clos span}_{\zeta\in\Pi^N}\{
[\varphi(\zeta)-\varphi(e)]\mathcal{U}\}
\perp\mbox{clos}\{\varphi(e)\mathcal{U}\} =:\mathcal{X}.$$ For any
$\zeta\in\Pi^N$ and $u\in\mathcal{U}$ one can represent now
$\varphi(\zeta)u$ as $$\mbox{col}\left[\begin{array}{cc}
  \varphi(e) & \varphi(\zeta)-\varphi(e)
\end{array}\right]u\in\mathcal{X\oplus H}.$$ On the other hand, for
any $u\in\mathcal{U},\ \zeta\in\Pi^N$ one has
$$\varphi(e)u\in\mbox{clos span}_
{\zeta\in\Pi^N}\{\varphi(\zeta )\mathcal{U}\},\quad [\varphi(\zeta
)-\varphi(e)]u\in\mbox{clos span}_ {\zeta\in\Pi^N}\{\varphi(\zeta
)\mathcal{U}\}.$$ Thus, $\mbox{clos
span}_{\zeta\in\Pi^N}\{\varphi(\zeta)\mathcal{U}\}=\mathcal{X\oplus
H}$. Let $\kappa:\mathcal{X\oplus H}\longrightarrow\mathcal{M}$ be
the natural embedding defined by
\begin{equation}\label{eq:kappa}
   \kappa:\left[\begin{array}{c}
  \varphi(e)u \\
(\varphi(\zeta )-\varphi(e))u
\end{array}\right]\longmapsto\varphi(\zeta )u=\left[\begin{array}{c}
  \varphi_1(\zeta )u \\
  \vdots \\
  \varphi_N(\zeta )u
\end{array}\right],\quad \zeta\in\Pi^N,\ u\in\mathcal{U},
\end{equation}
and extended to the whole $\mathcal{X\oplus H}$ by linearity and
continuity. Set
$$A_k:=\left[\begin{array}{cc}
  \varphi(e)^* & 0 \\
  0 & I_\mathcal{H}
\end{array}\right]\kappa^*P_k\kappa\left[\begin{array}{cc}
  \varphi(e) & 0 \\
  0 & I_\mathcal{H}
\end{array}\right]\in L\mathcal{(U\oplus H)},\quad k=1,\ldots,N.$$
Clearly, $A_k=A_k^*$ are PSD operators. Set
$$\psi(\zeta):=\left[\begin{array}{c}
  I_\mathcal{U}\\
\varphi(\zeta )-\varphi(e)
\end{array}\right]\in L\mathcal{(U,U\oplus H)},\quad
\zeta\in\Pi^N.$$ Then $f(z)=\psi(\zeta)^*A(z)\psi(z),\
(z,\zeta)\in\Pi^N\times\Pi^N$. Indeed,
\begin{eqnarray*}
\psi(\zeta)^*A(z)\psi(z)= \hspace{95mm}\\
\left[\begin{array}{c}
  I_\mathcal{U} \\
\varphi(\zeta)-\varphi(e)
\end{array}\right]^*\left[\begin{array}{cc}
  \varphi(e) & 0 \\
  0 & I_\mathcal{H}
\end{array}\right]^*\kappa^*\left(\sum_{k=1}^Nz_kP_k\right)
\kappa\left[\begin{array}{cc}
  \varphi(e) & 0 \\
  0 & I_\mathcal{H}
\end{array}\right]\left[\begin{array}{c}
  I_\mathcal{U} \\
\varphi(z)-\varphi(e)
\end{array}\right] \\
=\sum_{k=1}^Nz_k\varphi_k(\zeta)^*\varphi_k(z)=f(z).
\end{eqnarray*}
Now, let $A_k,\ k=1,\ldots,N$, have the block partitioning:
$$A_k=\left[\begin{array}{cc}
  a_k & b_k \\
  c_k & d_k
\end{array}\right]\in L\mathcal{(U\oplus H)}.$$
Then
\begin{eqnarray*}
A(z)\psi(z)=\left[\begin{array}{cc}
  a(z) & b(z) \\
  c(z) & d(z)
\end{array}\right]\left[\begin{array}{c}
  I_\mathcal{U} \\
\varphi(z)-\varphi(e)
\end{array}\right]=\\
\left[\begin{array}{c}
  a(z)+b(z)(\varphi(z)-\varphi(e)) \\
  c(z)+d(z)(\varphi(z)-\varphi(e))
\end{array}\right]=:\left[\begin{array}{c}
  f_1(z) \\
  f_2(z)
\end{array}\right].
\end{eqnarray*}
Since for any $(z,\zeta)\in\Pi^N\times\Pi^N$ one has
$$\psi(\zeta)^*A(z)\psi(z)=\left[\begin{array}{cc}
  I_\mathcal{U} & \varphi(\zeta)^*-\varphi(e)^*
\end{array}\right]\left[\begin{array}{c}
  f_1(z) \\
  f_2(z)
\end{array}\right]=f(z),$$
by setting $\zeta:=e$ in this equality we get $f_1(z)=f(z),\
z\in\Pi^N$. Therefore, for any $(z,\zeta)\in\Pi^N\times\Pi^N$ one
has $[\varphi(\zeta)-\varphi(e)]^*f_2(z)=0$. This implies that for
any $z\in\Pi^N$ and $u\in\mathcal{U}$ one has
$f_2(z)u\perp\mathcal{H}$. But $f_2(z)u\in\mathcal{H}$. Therefore,
$f_2(z)u=0$, and $f_2(z)\equiv 0$, i.e.,
\begin{equation}\label{eq:fi}
c(z)+d(z)[\varphi(z)-\varphi(e)]\equiv 0.
\end{equation}
Since for any $z\in\Pi^N$ the operator $P(z):=\sum_{k=1}^Nz_kP_k$
has positive definite real part, the operator
$d(z)=P_\mathcal{H}\kappa^*P(z)\kappa |\mathcal{H}$ has positive
definite real part, too. Therefore, $d(z)$ is boundedly invertible
for all $z\in\Pi^N$. From (\ref{eq:fi}) we get
$\varphi(z)-\varphi(e)=-d(z)^{-1}c(z),\ z\in\Pi^N$, and
$$f(z)=f_1(z)=a(z)-b(z)d(z)^{-1}c(z),\quad z\in\Pi^N.$$
Taking into account Remark~\ref{rem:dom}, we get
$f\in\mathcal{B}_N(\mathcal{U})$.

Since $\varphi(z)-\varphi(e)=-d(z)^{-1}c(z)$ and, hence, $\psi(z)$
are well-defined, holomorphic and homogeneous of degree zero
functions on $\Omega_N$, the kernels
$\Phi_k(z,\zeta)=\psi(\zeta)^*A_k\psi(z),\ k=1,\ldots, N$, are
extended to $\Omega_N\times\Omega_N$, and (\ref{eq:homzero})
holds. One can easily verify that these extended functions are
holomorphic PSD kernels on $\Omega_N\times\Omega_N$, and
(\ref{eq:b}) holds. The proof is complete.
\end{proof}

\section{The class $\mathcal{B}_N(\mathcal{U})$ and functional
calculus}\label{sec:b-calc} In this section we will give a
characterization of the class $\mathcal{B}_N(\mathcal{U})$ via the
functional calculus of $N$-tuples of commuting bounded strictly
accretive operators. First of all, let us observe that the
identity (\ref{eq:kerid}) is equivalent to the pair of the
following identities:
\begin{eqnarray}
f(z)+f(\zeta)^*=\sum\limits_{k=1}^N(z_k+\overline{\zeta_k})\Phi_k(z,\zeta),\quad
(z,\zeta)\in\Pi^N\times\Pi^N, \label{eq:b+}\\
f(z)-f(\zeta)^*=\sum\limits_{k=1}^N(z_k-\overline{\zeta_k})\Phi_k(z,\zeta),\quad
(z,\zeta)\in\Pi^N\times\Pi^N. \label{eq:b-}
\end{eqnarray}
We will show that the Cayley transform over all the variables
turns the first of these identities into Agler's identity which
characterizes the Agler--Herglotz class of holomorphic functions
on the unit polydisk $\mathbb{D}^N$. The latter has also, due to
\cite{Ag}, another characterization, via the functional calculus
of $N$-tuples of commuting strict contractions, that will give us
the desired result for the class $\mathcal{B}_N(\mathcal{U})$.

Let us recall the necessary definitions. Denote by $\mathcal{C}^N$
the set of all $N$-tuples $\mathbf{T}=(T_1,\ldots,T_N)$ of
commuting linear operators on some common Hilbert space
$\mathcal{H}$ subject to the condition $\| T_k\| <1,\
k=1,\ldots,N$ (\emph{strict contractions}). Then for any
holomorphic $L\mathcal{(U)}$-valued function
$$F(w)=\sum_{t\in\mathbb{Z}^N_+}\widehat{F_t}w^t,\quad
w\in\mathbb{D}^N$$ (here $\mathbb{Z}^N_+:=\{ t\in\mathbb{Z}^N:\
t_k\geq 0,k=1,\ldots,N\}$ and $w^t:=\prod_{k=1}^Nw_k^{t_k}$), and
any $\mathbf{T}\in\mathcal{C}^N$ the operator
\begin{equation}\label{eq:contr-calc}
    F(\mathbf{T}):=\sum_{t\in\mathbb{Z}^N_+}\widehat{F_t}\otimes
    \mathbf{T}^t\in L\mathcal{(U\otimes H)}
\end{equation}
is well-defined as a sum of a series convergent in the operator
norm. The \emph{Agler-Herglotz} class
$\mathcal{AH}_N(\mathcal{U})$ consists of all
$L\mathcal{(U)}$-valued functions $F$ which are holomorphic on
$\mathbb{D}^N$ and satisfying the condition
\begin{equation}\label{eq:AH}
    F(\mathbf{T})+F(\mathbf{T})^*\geq 0,\quad
    \mathbf{T}\in\mathcal{C}^N,
\end{equation}
where the inequality ``$\geq$" is considered in the sense of
positive semi-definiteness of a selfadjoint operator. It was
proved in \cite{Ag} that $F\in\mathcal{AH}_N(\mathcal{U})$ if and
only if there exist holomorphic PSD kernels $\Xi_k(w,\omega),\
k=1,\ldots,N$, on $\mathbb{D}^N\times\mathbb{D}^N$ such that
\begin{equation}\label{eq:AH-id}
    F(w)+F(\omega)^*=\sum_{k=1}^N(1-\overline{\omega_k}w_k)
    \Xi_k(w,\omega),\quad
    (w,\omega)\in\mathbb{D}^N\times\mathbb{D}^N.
\end{equation}
Denote by $\mathcal{A}^N$ the set of all $N$-tuples
$\mathbf{R}=(R_1,\ldots,R_N)$ of commuting bounded linear
operators on some common Hilbert space $\mathcal{H}$ for which
there exists a real constant $s>0$ such that $R_k+R_k^*\geq s
I_\mathcal{H},\ k=1,\ldots,N$ (\emph{strictly accretive
operators}). The \emph{operator Cayley transform}, defined by
\begin{equation}\label{eq:Cayley}
R_k:=(I_\mathcal{H}+T_k)(I_\mathcal{H}-T_k)^{-1},\quad
k=1,\ldots,N,
\end{equation}
maps the set $\mathcal{C}^N$ onto $\mathcal{A}^N$, and its inverse
\begin{equation}\label{eq:inv-Cayley}
T_k:=(R_k-I_\mathcal{H})(R_k+I_\mathcal{H})^{-1},\quad
k=1,\ldots,N,
\end{equation}
maps $\mathcal{A}^N$ onto $\mathcal{C}^N$.

Let $f$ be an $L\mathcal{(U)}$-valued function holomorphic on
$\Pi^N$. Then the \emph{Cayley transform over variables} applied
to $f$ gives
\begin{equation}\label{eq:F}
    F(w):=f\left(\frac{1+w_1}{1-w_1},\ldots,\frac{1+w_N}{1-w_N}\right),\quad
    w\in\mathbb{D}^N,
\end{equation}
which is a holomorphic $L\mathcal{(U)}$-valued function on
$\mathbb{D}^N$.

For any $\mathbf{R}\in\mathcal{A}^N$ let us define
$f(\mathbf{R}):=F(\mathbf{T})$ where $\mathbf{T}=(T_1,\ldots,T_N)$
is defined by (\ref{eq:inv-Cayley}).
\begin{thm}\label{thm:b-calc}
Let $f$ be an $L\mathcal{(U)}$-valued function holomorphic on
$\Pi^N$. Then $f\in\mathcal{B}_N(\mathcal{U})$ if and only if the
following conditions are satisfied:
\begin{description}
  \item[(i)] $f(tz_1,\ldots, tz_N)=tf(z_1,\ldots, z_N),\quad t>0,\
    z=(z_1,\ldots, z_N)\in\Pi^N$;
    \item[(ii)] $f(\mathbf{R})+f(\mathbf{R})^*\geq 0,\quad
    \mathbf{R}\in\mathcal{A}^N$;
    \item[(iii)] $f(\bar{z})=f(z)^*,\quad z\in\Pi^N$.
\end{description}
\end{thm}
\begin{proof}
 \textbf{Necessity.} Let $f\in\mathcal{B}_N(\mathcal{U})$. Then
(i) and (iii) easily follow from the representation (\ref{eq:gsc})
of $f$. Condition (ii) on $f$ is equivalent to condition
(\ref{eq:AH}) on $F$ which is defined by (\ref{eq:F}), i.e., to
$F\in\mathcal{AH}_N(\mathcal{U})$. Let us show the latter. Since
$f$ satisfies (\ref{eq:b+}), one can set
$z_k:=\frac{1+w_k}{1-w_k},\
\zeta_k=\frac{1+\omega_k}{1-\omega_k},\ k=1,\ldots,N$, in
(\ref{eq:b+}) and get
\begin{eqnarray}
\lefteqn{F(w)+F(\omega)^*=}\nonumber\\
& & \sum_{k=1}^N\left(\frac{1+w_k}{1-w_k}+
\frac{1+\overline{\omega_k}}{1-\overline{\omega_k}}\right)
\Phi_k\left(\frac{1+w_1}{1-w_1},\ldots,\frac{1+w_N}{1-w_N};
\frac{1+\omega_1}{1-\omega_1},\ldots,\frac{1+\omega_N}{1-\omega_N}\right)
\nonumber\\
& & \ \ \ \
=\sum_{k=1}^N(1-\overline{\omega_k}w_k)\Xi_k(w,\omega),\quad
(w,\omega)\in\mathbb{D}^N\times\mathbb{D}^N,\label{eq:AH+}
\end{eqnarray}
where for $k=1,\ldots,N$:
\begin{equation}\label{eq:AH-ker}
    \Xi_k(w,\omega)=\frac{2}{(1-w_k)(1-\overline{\omega_k})}
    \Phi_k\left(\frac{1+w_1}{1-w_1},\ldots,\frac{1+w_N}{1-w_N};
\frac{1+\omega_1}{1-\omega_1},\ldots,\frac{1+\omega_N}{1-\omega_N}\right).
\end{equation}
Since for each $k=1,\ldots,N$ one has
$\Phi_k(z,\zeta)=\varphi_k(\zeta)^*\varphi_k(z),\
(z,\zeta)\in\Pi^N\times\Pi^N$, where $\varphi_k$ is a holomorphic
function on $\Pi^N$ with values in $L\mathcal{(U,M}_k)$ for an
auxiliary Hilbert space $\mathcal{M}_k$ (again, see \cite{Aro}),
we get $\Xi_k(w,\omega)=\xi_k(\omega)^*\xi_k(w),\
(w,\omega)\in\mathbb{D}^N\times\mathbb{D}^N$, where
\begin{equation}\label{eq:xi}
    \xi_k(w)=\frac{\sqrt{2}}{1-w_k}\varphi_k
    \left(\frac{1+w_1}{1-w_1},\ldots,\frac{1+w_N}{1-w_N}\right),
\end{equation}
are holomorphic $L\mathcal{(U,M}_k)$-valued functions on
$\mathbb{D}^N$. Thus, $\Xi_k(w,\omega),\  k=1,\ldots,N$, are
    holomorphic PSD kernels on $\mathbb{D}^N\times\mathbb{D}^N$,
    and (\ref{eq:AH+}) means that $F\in\mathcal{AH}_N(\mathcal{U})$.

\textbf{Sufficiency.} Let $f$ satisfy conditions (i)--(iii). Since
(ii) is equivalent to $F\in\mathcal{AH}_N(\mathcal{U})$, where $F$
is defined by (\ref{eq:F}), the identity (\ref{eq:AH-id}) holds
with holomorphic PSD kernels $\Xi_k(w,\omega)$ on
$\mathbb{D}^N\times\mathbb{D}^N$. Let
$\Xi_k(w,\omega)=\xi_k(\omega)^*\xi_k(w),\
(w,\omega)\in\mathbb{D}^N\times\mathbb{D}^N$, where $\xi_k$ are
holomorphic functions on $\mathbb{D}^N$ taking values in
$L\mathcal{(U,M}_k)$ for some auxiliary Hilbert spaces
$\mathcal{M}_k,\ k=1,\ldots,N$. Set $w_k:=\frac{z_k-1}{z_k+1},\
\omega_k=\frac{\xi_k-1}{\xi_k+1},\ k=1,\ldots,N$, in
(\ref{eq:AH-id}), and by virtue of (\ref{eq:F}) get:
\begin{eqnarray}
\lefteqn{f(z)+f(\zeta)^*=}\nonumber\\
& &
\sum_{k=1}^N\left(1-\frac{\overline{\zeta_k}-1}{\overline{\zeta_k}+1}\cdot
\frac{z_k-1}{z_k+1}\right)
\xi_k\left(\frac{\zeta_1-1}{\zeta_1+1},\ldots,\frac{\zeta_N-1}{\zeta_N+1}\right)^*
\xi_k\left(\frac{z_1-1}{z_1+1},\ldots,\frac{z_N-1}{z_N+1}\right)
\nonumber\\
& &\ \ \ \ \
=\sum_{k=1}^N(z_k+\overline{\zeta_k})\varphi_k(\zeta)^*\varphi_k(z),\quad
(z,\zeta)\in\Pi^N\times\Pi^N,\label{eq:b++}
\end{eqnarray}
where for $k=1,\ldots,N$
\begin{equation}\label{eq:phi}
    \varphi_k(z):=\frac{\sqrt{2}}{z_k+1}
    \xi_k\left(\frac{z_1-1}{z_1+1},\ldots,\frac{z_N-1}{z_N+1}\right)
\end{equation}
are holomorphic $L\mathcal{(U,M}_k)$-valued functions on $\Pi^N$.
It follows that $\Phi_k(z,\zeta):=\varphi_k(\zeta)^*\varphi_k(z),\
k=1,\ldots,N$, are holomorphic PSD kernels on $\Pi^N\times\Pi^N$,
and (\ref{eq:b+})
  holds. The property (iii) implies $f(x)=f(x)^*$ for
any $x\in\mathbb{R}^N\cap\Pi^N$, and for any such $x$, and $t>0$
by (\ref{eq:b++}) one has:
\begin{eqnarray*}
f(x)+f(tx) &=& (1+t)\sum_{k=1}^Nx_k\varphi_k(tx)^*\varphi_k(x),\\
f(tx)+f(x) &=& (1+t)\sum_{k=1}^Nx_k\varphi_k(x)^*\varphi_k(tx),\\
\frac{1+t}{2}[f(x)+f(x)] &=& \frac{1+t}{2}\sum_{k=1}^N2x_k
\varphi_k(x)^*\varphi_k(x),\\
\frac{1+t}{2t}[f(tx)+f(tx)] &=& \frac{1+t}{2t}
\sum_{k=1}^N2tx_k\varphi_k(tx)^*\varphi_k(tx).
\end{eqnarray*}
By (i), the left-hand sides of these equalities coincide and equal
to $(1+t)f(x)$, therefore
\begin{eqnarray*}
f(x) &=& \sum_{k=1}^Nx_k\varphi_k(tx)^*\varphi_k(x)=
\sum_{k=1}^Nx_k\varphi_k(x)^*\varphi_k(tx)\\
     &=& \sum_{k=1}^Nx_k
\varphi_k(x)^*\varphi_k(x)=\sum_{k=1}^Nx_k\varphi_k(tx)^*\varphi_k(tx).
\end{eqnarray*}
>From here we get
\begin{eqnarray*}
0 &\leq &
\sum_{k=1}^Nx_k[\varphi_k(tx)-\varphi_k(x)]^*[\varphi_k(tx)-\varphi_k(x)]\\
  &=& \sum_{k=1}^Nx_k\varphi_k(tx)^*\varphi_k(tx)-
  \sum_{k=1}^Nx_k\varphi_k(tx)^*\varphi_k(x)\\
  &-& \sum_{k=1}^Nx_k\varphi_k(x)^*\varphi_k(tx)+
  \sum_{k=1}^Nx_k\varphi_k(x)^*\varphi_k(x)=0.
\end{eqnarray*}
Thus $\varphi_k(tx)-\varphi_k(x)=0$ for any
$x\in\mathbb{R}^N\cap\Pi^N,\ t>0$ and $k=1,\ldots,N$. For fixed
$k\in\{ 1,\ldots,N\}$ and $t>0$ the function
$h_{k,t}(z):=\varphi_k(tz)-\varphi_k(z)$ is holomorphic on $\Pi^N$
and takes values in $L\mathcal{(U,M}_k)$. Then for any fixed
$k\in\{ 1,\ldots,N\},\ t>0,\ u\in\mathcal{U}$ and
$m\in\mathcal{M}_k$ the scalar function $h_{k,t,u,m}(z):=\langle
h_{k,t}(z)u,m\rangle_{\mathcal{M}_k}$ is holomorphic on $\Pi^N$
and vanishes on $\mathbb{R}^N\cap\Pi^N$. By the uniqueness theorem
for holomorphic functions of several variables (see, e.g.,
\cite{Sh2}), $h_{k,t,u,m}(z)\equiv 0$, hence $h_{k,t}(z)\equiv 0$,
that means:
$$\varphi_k(tz)=\varphi_k(z),\quad t>0,\ z\in\Pi^N.$$
It follows from here that for any $(z,\zeta)\in\Pi^N\times\Pi^N$
and $t>0$ one has
\begin{eqnarray*}
f(z)+tf(\zeta)^* &=& f(z)+f(t\zeta)^*=
\sum_{k=1}^N(z_k+t\overline{\zeta_k})\varphi_k(t\zeta)^*\varphi_k(z)\\
 &=&
 \sum_{k=1}^N(z_k+t\overline{\zeta_k})\varphi_k(\zeta)^*\varphi_k(z)\\
 &=&
\sum_{k=1}^Nz_k\varphi_k(\zeta)^*\varphi_k(z)+
t\sum_{k=1}^N\overline{\zeta_k}\varphi_k(\zeta)^*\varphi_k(z),
\end{eqnarray*}
and the comparison of the coefficients of the two linear functions
in $t$, in the beginning and in the end of this chain of
equalities, gives:
$$f(z)=\sum_{k=1}^Nz_k\varphi_k(\zeta)^*\varphi_k(z),\quad
(z,\zeta)\in\Pi^N\times\Pi^N,$$ i.e., (\ref{eq:kerid}). By
Theorem~\ref{thm:ker}, $f\in\mathcal{B}_N(\mathcal{U})$. The proof
is complete.
\end{proof}
\begin{cor}\label{cor:bes}
Let $f$ be an $L\mathcal{(U)}$-valued function holomorphic on
$\Omega_N$. Then $f\in\mathcal{B}_N(\mathcal{U})$ if and only if
the following conditions are satisfied:
\begin{description}
  \item[(i)] $f(\lambda z_1,\ldots, \lambda z_N)=
  \lambda f(z_1,\ldots, z_N),\quad \lambda\in\mathbb{C}\backslash\{ 0\},\
    z=(z_1,\ldots, z_N)\in\Omega_N$;
    \item[(ii)] $f(\mathbf{R})+f(\mathbf{R})^*\geq 0,\quad
    \mathbf{R}\in\mathcal{A}^N$;
    \item[(iii)] $f(\bar{z})=f(z)^*,\quad z\in\Omega_N$.
\end{description}
\end{cor}
\begin{proof}
If $f\in\mathcal{B}_N(\mathcal{U})$ then conditions (i) and (iii)
follow from the representation (\ref{eq:gsc}) of $f$, and
condition (ii) follows from Theorem~\ref{thm:b-calc}. Conversely,
conditions (i)--(iii) of this Corollary imply conditions
(i)--(iii) of Theorem~\ref{thm:b-calc}, which in turn imply
$f\in\mathcal{B}_N(\mathcal{U})$.
\end{proof}
\begin{rem}\label{rem:on(ii)}
By Corollary~\ref{cor:bes}, its conditions (i)--(iii) on
holomorphic $L\mathcal{(U)}$-valued functions on $\Omega_N$ give
an equivalent definition of the class
$\mathcal{B}_N(\mathcal{U})$, which seems to be more natural than
the original definition given in Section~\ref{sec:classes} in
``existence" terms. The definition of the class
$\mathcal{P}_N(\mathcal{U})$ is obtained by replacing condition
(ii) by a weaker condition (\ref{eq:gpos}).
\end{rem}

\section{The image of the Bessmertny\u{\i} class under the double Cayley
transform}\label{sec:Ag} In Section~\ref{sec:b-calc} it was
established that the Cayley transform over the variables maps the
Bessmertny\u{\i} class $\mathcal{B}_N(\mathcal{U})$ into the
Agler--Herglotz class  $\mathcal{AH}_N(\mathcal{U})$. Since the
Cayley transform over the values of functions maps the class
$\mathcal{AH}_N(\mathcal{U})$ into the Agler--Schur class
$\mathcal{AS}_N(\mathcal{U})$ (see \cite{Ag}), the composition of
these two transforms (the \emph{double Cayley transform}) maps
$\mathcal{B}_N(\mathcal{U})$ into $\mathcal{AS}_N(\mathcal{U})$.
The class $\mathcal{AS}_N(\mathcal{U})$ is important in the
interpolation theory and systems theory in several variables (see,
e.g., \cite{Ag, AgMC, BT, K2, K3, BLTT}), that is why it is
interesting to describe the image of $\mathcal{B}_N(\mathcal{U})$
in $\mathcal{AS}_N(\mathcal{U})$ under the double Cayley
transform.

Given $f\in\mathcal{B}_N(\mathcal{U})$, define for
$w\in\mathbb{D}^N$:
\begin{eqnarray}
\lefteqn{\mathcal{F}(w)=(F(w)-I_\mathcal{U})(F(w)+I_\mathcal{U})^{-1}=}\nonumber\\
& & \left(
f\left(\frac{1+w_1}{1-w_1},\ldots,\frac{1+w_N}{1-w_N}\right)
-I_\mathcal{U}\right)\left(
f\left(\frac{1+w_1}{1-w_1},\ldots,\frac{1+w_N}{1-w_N}\right)
+I_\mathcal{U}\right)^{-1}. \label{eq:Cayley-val}
\end{eqnarray}
We shall write down $\mathcal{F}=\mathcal{C}(f)$, and call
$\mathcal{C}(\cdot)$ the \emph{double Cayley transform}.

Let us recall the definition, and  resume the main results of
\cite{Ag} on the Agler--Schur class $\mathcal{AS}_N(\mathcal{U})$.
\begin{thm}\label{thm:Ag}
Let $\mathcal{F}$ be a holomorphic $L\mathcal{(U)}$-valued
function on $\mathbb{D}^N$. The following statements are
equivalent:
\begin{description}
    \item[(i)] $\|\mathcal{F}(\mathbf{T})\|\leq 1$ for any
    $\mathbf{T}\in\mathcal{C}^N$;
    \item[(ii)] there exist holomorphic PSD kernels
    $\Theta_k(w,\omega)$ on $\mathbb{D}^N\times\mathbb{D}^N,\
    k=1,\ldots,N$, such that for any
    $(w,\omega)\in\mathbb{D}^N\times\mathbb{D}^N$ one has
    \begin{equation}\label{eq:Ag-id}
    I_\mathcal{U}-\mathcal{F(\omega})^*\mathcal{F}(w)=
    \sum_{k=1}^N(1-\overline{\omega_k}w_k)\Theta_k(w,\omega);
\end{equation}
    \item[(iii)] there exist Hilbert spaces
    $\mathcal{X,X}_1,\ldots,\mathcal{X}_N$ such that
    $\mathcal{X}=\bigoplus_{k=1}^N\mathcal{X}_k$, and a unitary
    operator $U=\left[\begin{array}{cc}
      A & B \\
      C & D
    \end{array}\right]\in L\mathcal{(X\oplus U)}$ such that
    \begin{equation}\label{eq:Ag}
    \mathcal{F}(w)=D+CP(w)(I_\mathcal{X}-AP(w))^{-1}B,
\end{equation}
where $P(w):=\sum_{k=1}^Nw_kP_{\mathcal{X}_k}$, i.e.,
$\mathcal{F}$ is the \textbf{transfer function} of an
\textbf{Agler unitary colligation}
$\alpha=(N;U;\mathcal{X}=\bigoplus_{k=1}^N\mathcal{X}_k,\mathcal{U},\mathcal{U})$
(we will write $\mathcal{F=F}_\alpha$ in this case) with the
\textbf{state space} $\mathcal{X}$, and the same \textbf{input and
output spaces} equal to $\mathcal{U}$.
\end{description}
\end{thm}
The \emph{Agler--Schur class} $\mathcal{AS}_N(\mathcal{U})$
consists of all functions satisfying any (and, hence, all) of
conditions (i)--(iii) of Theorem~\ref{thm:Ag}.
\begin{thm}\label{thm:double-C}
A holomorphic $L\mathcal{(U)}$-valued function $\mathcal{F}$ on
$\mathbb{D}^N$ can be represented as $\mathcal{F}=\mathcal{C}(f)$
for some $f\in\mathcal{B}_N(\mathcal{U})$ if and only if the
following conditions are fulfilled:
\begin{description}
    \item[(i)] $\mathcal{F}=\mathcal{F}_\alpha$ for an Agler unitary colligation
$\alpha=(N;U;\mathcal{X}=\bigoplus_{k=1}^N\mathcal{X}_k,\mathcal{U},\mathcal{U})$
with the additional condition $U=U^*$;
    \item[(ii)] $1\notin\sigma(\mathcal{F}(0))$.
\end{description}
\end{thm}
\begin{proof}
\textbf{Necessity.} Let $f\in\mathcal{B}_N(\mathcal{U})$. Then
(\ref{eq:b+}) and (\ref{eq:b-}) hold. As we have shown in
Theorem~\ref{thm:b-calc}, the identity (\ref{eq:b+}) implies the
identity (\ref{eq:AH+}) for the holomorphic
$L\mathcal{(U)}$-valued function $F$ on $\mathbb{D}^N$ which is
defined by (\ref{eq:F}), with holomorphic PSD kernels
$\Xi_k(w,\omega),\ k=1,\ldots,N$, on
$\mathbb{D}^N\times\mathbb{D}^N$ defined by (\ref{eq:AH-ker}).
Analogously, the identity (\ref{eq:b-}) implies
\begin{eqnarray}
\lefteqn{F(w)-F(\omega)^*=}\nonumber\\
& & \sum_{k=1}^N\left(\frac{1+w_k}{1-w_k}-
\frac{1+\overline{\omega_k}}{1-\overline{\omega_k}}\right)
\Phi_k\left(\frac{1+w_1}{1-w_1},\ldots,\frac{1+w_N}{1-w_N};
\frac{1+\omega_1}{1-\omega_1},\ldots,\frac{1+\omega_N}{1-\omega_N}\right)
\nonumber\\
& & \ \ \ \
=\sum_{k=1}^N(w_k-\overline{\omega_k})\Xi_k(w,\omega),\quad
(w,\omega)\in\mathbb{D}^N\times\mathbb{D}^N,\label{eq:AH-}
\end{eqnarray}
with the same kernels $\Xi_k(w,\omega),\ k=1,\ldots,N$. Let
$\mathcal{F}=\mathcal{C}(f)$, i.e., $\mathcal{F}$ is determined by
$f$ or $F$ in accordance with (\ref{eq:Cayley-val}). Then
\begin{eqnarray*}
I_\mathcal{U}-\mathcal{F}(\omega)^*\mathcal{F}(w)=
I_\mathcal{U}-\left(F(\omega)^*+I_\mathcal{U}\right)^{-1}\left(F(\omega)^*-
I_\mathcal{U}\right)
\left(F(w)-I_\mathcal{U}\right)\left(F(w)+I_\mathcal{U}\right)^{-1}\\
=\left(F(\omega)^*+I_\mathcal{U}\right)^{-1}\left[
\left(F(\omega)^*+I_\mathcal{U}\right)\left(F(w)+I_\mathcal{U}\right)
-\left(F(\omega)^*-I_\mathcal{U}\right)\left(F(w)-I_\mathcal{U}\right)\right]\\
\times\left(F(w)+I_\mathcal{U}\right)^{-1}
=2\left(F(\omega)^*+I_\mathcal{U}\right)^{-1}\left(F(w)+F(\omega)^*\right)
\left(F(w)+I_\mathcal{U}\right)^{-1}.
\end{eqnarray*}
According to (\ref{eq:AH+}), we get
\begin{equation}\label{eq:Ag+}
    I_\mathcal{U}-\mathcal{F}(\omega)^*\mathcal{F}(w)=
    \sum_{k=1}^N(1-\overline{\omega_k}w_k)\Theta_k(w,\omega),\quad
    (w,\omega)\in\mathbb{D}^N\times\mathbb{D}^N,
\end{equation}
where
\begin{equation}\label{eq:ker-Ag}
    \Theta_k(w,\omega)=2\left(F(\omega)^*+I_\mathcal{U}\right)^{-1}
    \Xi_k(w,\omega)\left(F(w)+I_\mathcal{U}\right)^{-1}.
\end{equation}
Since $\Xi_k(w,\omega)=\xi_k(\omega)^*\xi_k(w),\
    (w,\omega)\in\mathbb{D}^N\times\mathbb{D}^N$, where $\xi_k$
    are holomorphic $L\mathcal{(U,M}_k)$-valued functions on
    $\mathbb{D}^N$, one has
\begin{equation}\label{eq:fact-Ag}
    \Theta_k(w,\omega)=\theta_k(\omega)^*\theta_k(w),\quad
    (w,\omega)\in\mathbb{D}^N\times\mathbb{D}^N,
\end{equation}
where
\begin{equation}\label{eq:theta}
    \theta_k(w)=\sqrt{2}\xi_k(w)\left(F(w)+I_\mathcal{U}\right)^{-1},
\quad (w,\omega)\in\mathbb{D}^N\times\mathbb{D}^N,
\end{equation}
are also holomorphic $L\mathcal{(U,M}_k)$-valued functions on
    $\mathbb{D}^N$. Thus, $\Theta_k(w,\omega),\ k=1,\ldots,N$, are
    holomorphic PSD kernels on $\mathbb{D}^N\times\mathbb{D}^N$,
    and due to (\ref{eq:Ag+}),
    $\mathcal{F}\in\mathcal{AS}_N(\mathcal{U})$. Analogously,
    $$\mathcal{F}(w)-\mathcal{F}(\omega)^*
    =2\left(F(\omega)^*+I_\mathcal{U}\right)^{-1}\left(F(w)-F(\omega)^*\right)
\left(F(w)+I_\mathcal{U}\right)^{-1},$$ and according to
(\ref{eq:AH-}) we get
\begin{equation}\label{eq:Ag-}
   \mathcal{F}(w)- \mathcal{F}(\omega)^*=
    \sum_{k=1}^N(w_k-\overline{\omega_k})\Theta_k(w,\omega),\quad
    (w,\omega)\in\mathbb{D}^N\times\mathbb{D}^N,
\end{equation}
with the same set of kernels $\Theta_k(w,\omega),\ k=1,\ldots,N$,
defined by (\ref{eq:ker-Ag}). Let us rewrite (\ref{eq:Ag+}) and
(\ref{eq:Ag-}) in a somewhat different way. Since
$f\in\mathcal{B}_N(\mathcal{U})$ satisfies
$f(\overline{z})=f(z)^*,\ z\in\Pi^N$, one has also
$$F(\overline{w})=F(w)^*,\quad
\mathcal{F}(\overline{w})=\mathcal{F}(w)^*,\quad
w\in\mathbb{D}^N.$$ Therefore, (\ref{eq:Ag+}) and (\ref{eq:Ag-})
are equivalent to the following two identities, respectively:
\begin{equation}\label{eq:Ag++}
I_\mathcal{U}-\mathcal{F}(w)\mathcal{F}(\omega)^*=
    \sum_{k=1}^N(1-w_k\overline{\omega_k})
    \widetilde{\theta_k}(w)\widetilde{\theta_k}(\omega)^*,\quad
    (w,\omega)\in\mathbb{D}^N\times\mathbb{D}^N,
\end{equation}
\begin{equation}\label{eq:Ag--}
\mathcal{F}(w)-\mathcal{F}(\omega)^*=
    \sum_{k=1}^N(w_k-\overline{\omega_k})
    \widetilde{\theta_k}(w)\widetilde{\theta_k}(\omega)^*,\quad
    (w,\omega)\in\mathbb{D}^N\times\mathbb{D}^N,
\end{equation}
where $\widetilde{\theta_k}(w)=\theta_k(\bar{w})^*$ are
holomorphic $L\mathcal{(M}_k,\mathcal{U})$-valued functions on
$\mathbb{D}^N$. We will show that the identities (\ref{eq:Ag++})
and (\ref{eq:Ag--}) allow to construct an Agler unitary
colligation satisfying condition (i) of this theorem. For this
purpose, we will use the functional model of an Agler colligation
by J.~A.~Ball and T.~T.~Trent \cite{BT}.

Let us remind their construction. Let
$W\in\mathcal{AS}_N(\mathcal{U,Y})$ (the definition of
$\mathcal{AS}_N(\mathcal{U,Y})$ is the same as of
$\mathcal{AS}_N(\mathcal{U})$, with only difference that values of
functions from this class are in $L\mathcal{(U,Y)}$). Then there
exist Hilbert spaces $\mathcal{L}_k$ and holomorphic functions
$H_k,H_{*k}$ on $\mathbb{D}^N$ taking values in
$L(\mathcal{L}_k,\mathcal{Y})$ and $L(\mathcal{L}_k,\mathcal{U}),\
k=1,\ldots,N$, respectively, such that
\begin{eqnarray}
\lefteqn{\left[\begin{array}{cc}
  I_\mathcal{U}-W(\bar{w})^*W(\bar{\omega}) & W(\bar{w})^*-W(\omega)^* \\
  W(w)-W(\bar{\omega}) &
  I_\mathcal{Y}-W(w)W(\omega)^*
\end{array}\right]=}\nonumber\\
& & \sum_{k=1}^N\left[\begin{array}{cc}
  1-w_k\overline{\omega_k} & w_k-\overline{\omega_k} \\
  w_k-\overline{\omega_k} & 1-w_k\overline{\omega_k}
\end{array}\right]\circ\left[\begin{array}{cc}
  H_{*k}(w)H_{*k}(\omega)^* & H_{*k}(w)H_k(\omega)^* \\
  H_k(w)H_{*k}(\omega)^* & H_k(w)H_k(\omega)^*
\end{array}\right] \label{eq:Ag+-}
\end{eqnarray}
for all $(w,\omega)\in\mathbb{D}^N\times\mathbb{D}^N$, where
``$\circ$" is a \emph{Schur (entry-wise) matrix multiplication}.
For every $k=1,\ldots,N$ and
$(w,\omega)\in\mathbb{D}^N\times\mathbb{D}^N$, set
\begin{eqnarray*}
K_{*k}(w,\omega):=H_{*k}(w)H_{*k}(\omega)^*, &
K_k(w,\omega):=H_k(w)H_k(\omega)^*,\\
L_k(w,\omega):=H_{*k}(w)H_k(\omega)^*, &
L_{*k}(w,\omega):=L_k(\omega,w)^*=H_k(w)H_{*k}(\omega)^*,
\end{eqnarray*}
and $$ \widehat{K_k}(w,\omega):=\left[\begin{array}{cc}
  K_{*k}(w,\omega) & L_k(w,\omega) \\
  L_{*k}(w,\omega) & K_k(w,\omega)
\end{array}\right].$$
The latter function is a holomorphic $L(\mathcal{U\oplus
Y})$-valued PSD kernel on $\mathbb{D}^N\times\mathbb{D}^N$, which
serves as the reproducing kernel of the Hilbert space
$\mathcal{H}(\widehat{K_k})$ of holomorphic $(\mathcal{U\oplus
Y})$-valued functions on $\mathbb{D}^N$; this space is determined
uniquely by $\widehat{K_k}$ (for the theory of reproducing kernel
Hilbert spaces see, e.g., \cite{Aro}). Set
$\mathcal{D}(\{\widehat{K_k}\}_{k=1}^N):=
\bigoplus_{k=1}^N\mathcal{H}(\widehat{K_k})$. The latter is a
Hilbert space of holomorphic $(\mathcal{U\oplus Y})^N$-valued
functions on $\mathbb{D}^N$, with the reproducing kernel
$\widehat{K}(w,\omega):=\bigoplus_{k=1}^N\widehat{K_k}(w,\omega)$.
Define the lineals
$\mathcal{D}_0\subset\mathcal{D}(\{\widehat{K_k}\}_{k=1}^N)\oplus\mathcal{U}$
and
$\mathcal{R}_0\subset\mathcal{D}(\{\widehat{K_k}\}_{k=1}^N)\oplus\mathcal{Y}$
as
$$\mathcal{D}_0:=\mbox{span}\left\{\left[\begin{array}{c}
  \overline{w_1}K_{*1}(\cdot,w) \\
  \overline{w_1}L_{*1}(\cdot,w) \\
  \vdots \\
  \overline{w_N}K_{*N}(\cdot,w) \\
  \overline{w_N}L_{*N}(\cdot,w) \\
  I_\mathcal{U}
\end{array}\right]u,\ \left[\begin{array}{c}
  L_1(\cdot,w) \\
  K_1(\cdot,w) \\
  \vdots \\
  L_N(\cdot,w) \\
  K_N(\cdot,w) \\
  W(w)^*
\end{array}\right]y:\
u\in\mathcal{U},y\in\mathcal{Y}, w\in\mathbb{D}^N\right\},$$
$$\mathcal{R}_0:=\mbox{span}\left\{\left[\begin{array}{c}
  K_{*1}(\cdot,w) \\
  L_{*1}(\cdot,w) \\
  \vdots \\
  K_{*N}(\cdot,w) \\
  L_{*N}(\cdot,w) \\
  W(\overline{w})
\end{array}\right]u,\ \left[\begin{array}{c}
  \overline{w_1}L_1(\cdot,w) \\
  \overline{w_1}K_1(\cdot,w) \\
  \vdots \\
  \overline{w_N}L_N(\cdot,w) \\
  \overline{w_N}K_N(\cdot,w) \\
  I_\mathcal{Y}
\end{array}\right]y:\
u\in\mathcal{U},y\in\mathcal{Y}, w\in\mathbb{D}^N\right\}.$$ The
operator $U_0:\mathcal{D}_0\rightarrow\mathcal{R}_0$, correctly
defined by
$$U_0:\left[\begin{array}{c}
  \overline{w_1}K_{*1}(\cdot,w) \\
  \overline{w_1}L_{*1}(\cdot,w) \\
  \vdots \\
  \overline{w_N}K_{*N}(\cdot,w) \\
  \overline{w_N}L_{*N}(\cdot,w) \\
  I_\mathcal{U}
\end{array}\right]u\longmapsto\left[\begin{array}{c}
  K_{*1}(\cdot,w) \\
  L_{*1}(\cdot,w) \\
  \vdots \\
  K_{*N}(\cdot,w) \\
  L_{*N}(\cdot,w) \\
  W(\overline{w})
\end{array}\right]u,\quad u\in\mathcal{U}, w\in\mathbb{D}^N,$$
$$U_0:\left[\begin{array}{c}
  L_1(\cdot,w) \\
  K_1(\cdot,w) \\
  \vdots \\
  L_N(\cdot,w) \\
  K_N(\cdot,w) \\
  W(w)^*
\end{array}\right]y\longmapsto\left[\begin{array}{c}
  \overline{w_1}L_1(\cdot,w) \\
  \overline{w_1}K_1(\cdot,w) \\
  \vdots \\
  \overline{w_N}L_N(\cdot,w) \\
  \overline{w_N}K_N(\cdot,w) \\
  I_\mathcal{Y}
\end{array}\right]y,\quad
y\in\mathcal{Y}, w\in\mathbb{D}^N,$$ is uniquely extended to the
(correctly defined) unitary operator
$\widetilde{U_0}:\mbox{clos}(\mathcal{D}_0)\rightarrow\mbox{clos}(\mathcal{R}_0)$.
In the case where
\begin{equation}\label{eq:dim}
    \dim\{
    (\mathcal{D}(\{\widehat{K_k}\}_{k=1}^N)\oplus\mathcal{U})
    \ominus\mbox{clos}(\mathcal{D}_0)\}=\dim\{
    (\mathcal{D}(\{\widehat{K_k}\}_{k=1}^N)\oplus\mathcal{Y})
    \ominus\mbox{clos}(\mathcal{R}_0)\}
\end{equation}
there exists a (non-unique!) unitary operator
$U:\mathcal{D}(\{\widehat{K_k}\}_{k=1}^N)\oplus\mathcal{U}
\rightarrow\mathcal{D}(\{\widehat{K_k}\}_{k=1}^N)\oplus\mathcal{Y}$
such that $U|\mbox{clos}(\mathcal{D}_0)=\widetilde{U_0}$. The
corresponding Agler unitary colligation is
$\alpha=(N;U;\mathcal{D}(\{\widehat{K_k}\}_{k=1}^N)=
\bigoplus_{k=1}^N\mathcal{H}(\widehat{K_k}),U,Y)$, and
$W=W_\alpha$, i.e., $W$ is the transfer function of $\alpha$.

Let us apply this construction to a function
$\mathcal{F}=\mathcal{C}(f)$ where $f\in\mathcal{B}_N(U)$. In this
case $\mathcal{Y}=\mathcal{U},\ W=\mathcal{F}$. Since
$\mathcal{F}(\overline{w})=\mathcal{F}(w)^*,\ w\in\mathbb{D}^N$,
it follows from (\ref{eq:Ag++}) and (\ref{eq:Ag--}) that
(\ref{eq:Ag+-}) holds with $H_k=H_{*k}=\widetilde{\theta_k},\
k=1,\ldots,N$. Therefore, for $k=1,\ldots,N$ one has
\begin{equation}\label{eq:ker-eq}
    K_k(w,\omega)=K_{*k}(w,\omega)=L_k(w,\omega)=L_{*k}(w,\omega),\
    (w,\omega)\in\mathbb{D}^N\times\mathbb{D}^N.
\end{equation}
In turn, this means $\mathcal{D}_0=\mathcal{R}_0$, and $$
    (\mathcal{D}(\{\widehat{K_k}\}_{k=1}^N)\oplus\mathcal{U})
    \ominus\mbox{clos}(\mathcal{D}_0)=
    (\mathcal{D}(\{\widehat{K_k}\}_{k=1}^N)\oplus\mathcal{Y})
    \ominus\mbox{clos}(\mathcal{R}_0).$$
In particular, (\ref{eq:dim}) holds. Define the operator
$$U:=\widetilde{U_0}\oplus I_{(\mathcal{D}(\{\widehat{K_k}\}_{k=1}^N)
\oplus\mathcal{U})
    \ominus\mbox{clos}(\mathcal{D}_0)}\in L(\mathcal{D}(\{\widehat{K_k}\}_{k=1}^N)
\oplus\mathcal{U}).$$ Clearly, $U$ is unitary. Let us show that
$U=U^*$. From the definition of $U_0$ and (\ref{eq:ker-eq}) it
follows that $U_0=U_0^{-1}$. Therefore
$\widetilde{U_0}=\widetilde{U_0}^{-1}$ and $U=U^{-1}$. Since
$U^*=U^{-1}$, we obtain $U=U^*$. Thus, condition (i) is fulfilled.

Since $f(e)=f(e)^*\geq 0$ where $e=(1,\ldots,1)$,  the operator
$f(e)\in L(\mathcal{U})$ has a spectral decomposition (see, e.g.,
\cite{AkG})
$$f(e)=\int_0^{\| f(e)\|}t\,dE_t,$$
therefore,
$$\mathcal{F}(0)=(f(e)-I_\mathcal{U})(f(e)+I_\mathcal{U})^{-1}=
\int_0^{\| f(e)\|}\frac{t-1}{t+1}\,dE_t.$$ Since the function
$s(t)=\frac{t-1}{t+1}$ increases on the segment $[0,\| f(e)\| ]$,
one has $\mathcal{F}(0)\leq\frac{\| f(e)\| -1}{\| f(e)\|
+1}I_\mathcal{U}$. Since $\frac{\| f(e)\| -1}{\| f(e)\| +1}<1$, we
conclude that $1\notin\sigma(\mathcal{F}(0))$, i.e., condition
(ii) is also fulfilled.

\textbf{Sufficiency.} Let conditions (i) and (ii) for the function
$\mathcal{F}$ satisfy. Since $U=U^*=U^{-1}$, one has $D=D^*$ and
$D$ is a contraction in $\mathcal{U}$. Since $\mathcal{F}(0)=D$
and $1\notin\sigma(\mathcal{F}(0))$, we get
$\sigma(\mathcal{F}(0))\subset [-1,a]$ with some $a:\ -1\leq a<1$.
Hence, $\sigma(\frac{I_\mathcal{U}+\mathcal{F}(0)}{2})\subset
[0,\frac{1+a}{2}]$, thus
$\left\|\frac{I_\mathcal{U}+\mathcal{F}(0)}{2}\right\|\leq\frac{1+a}{2}<1$.
By the maximum principle for holomorphic operator-valued functions
of several variables (e.g., see \cite{Schw}),
$\left\|\frac{I_\mathcal{U}+\mathcal{F}(w)}{2}\right\| <1$ for all
$w\in\mathbb{D}^N$. Indeed, (i) implies
$\mathcal{F}\in\mathcal{AS}_N(\mathcal{U})$, and therefore
$\frac{I_\mathcal{U}+\mathcal{F}}{2}\in\mathcal{AS}_N(\mathcal{U})$.
Thus $\left\|\frac{I_\mathcal{U}+\mathcal{F}(w)}{2}\right\| \leq
1$ for all $w\in\mathbb{D}^N$. If for some $w_0\in\mathbb{D}^N$
one had $\left\|\frac{I_\mathcal{U}+\mathcal{F}(w_0)}{2}\right\|
=1$, then the maximum principle would imply
$\left\|\frac{I_\mathcal{U}+\mathcal{F}(w)}{2}\right\| =1$
everywhere in $\mathbb{D}^N$. In particular,
$\left\|\frac{I_\mathcal{U}+\mathcal{F}(0)}{2}\right\| =1$, that
is not true. Finally, we get
$1\notin\sigma(\frac{I_\mathcal{U}+\mathcal{F}(w)}{2})$, and
therefore $1\notin\sigma(\mathcal{F}(w))$ for all
$w\in\mathbb{D}^N$. Thus, the function
$F(w)=(I_\mathcal{U}+\mathcal{F}(w))(I_\mathcal{U}-\mathcal{F}(w))^{-1}$
is correctly defined and holomorphic on $\mathbb{D}^N$. It is easy
to see that
\begin{equation}\label{eq:F+}
    F(w)+F(\omega)^*=
    2(I_\mathcal{U}-\mathcal{F}(\omega)^*)^{-1}
    (I_\mathcal{U}-\mathcal{F}(\omega)^*\mathcal{F}(w))
    (I_\mathcal{U}-\mathcal{F}(w))^{-1},
\end{equation}
\begin{equation}\label{eq:F-}
    F(w)-F(\omega)^*=
    2(I_\mathcal{U}-\mathcal{F}(\omega)^*)^{-1}
    (\mathcal{F}(w)-\mathcal{F}(\omega)^*)
    (I_\mathcal{U}-\mathcal{F}(w))^{-1}
\end{equation}
for all $(w,\omega)\in\mathbb{D}^N\times\mathbb{D}^N$. Since
$U=U^{-1}=U^*$, due to (\ref{eq:Ag}) one has
\begin{eqnarray*}
\lefteqn{I_\mathcal{U}-\mathcal{F}(\omega)^*\mathcal{F}(w) =
I_\mathcal{U}-[D+CP(\omega)(I_\mathcal{X}-AP(\omega))^{-1}B]^*} \\
& & \times  [D+CP(w)(I_\mathcal{X}-AP(w))^{-1}B]=
I_\mathcal{U}-D^*D\\
& & -D^*C P(w)(I_\mathcal{X}-AP(w))^{-1}B
 - B^*(I_\mathcal{X}-P(\bar{\omega})A^*)^{-1}P(\bar{\omega})C^*D\\
 & &-B^*(I_\mathcal{X}-P(\bar{\omega})A^*)^{-1}P(\bar{\omega})C^*
 C P(w)(I_\mathcal{X}-AP(w))^{-1}B\\
 & & =B^*B+B^*AP(w)(I_\mathcal{X}-AP(w))^{-1}B+
 B^*(I_\mathcal{X}-P(\bar{\omega})A^*)^{-1}P(\bar{\omega})A^*B\\
 & & -B^*(I_\mathcal{X}-P(\bar{\omega})A^*)^{-1}P(\bar{\omega})
 (I_\mathcal{X}-A^*A)P(w)(I_\mathcal{X}-AP(w))^{-1}B\\
 & & =B^*(I_\mathcal{X}-P(\bar{\omega})A)^{-1}[(I_\mathcal{X}-P(\bar{\omega})A)
(I_\mathcal{X}-AP(w))+(I_\mathcal{X}-P(\bar{\omega})A)AP(w)\\
& &
+P(\bar{\omega})A(I_\mathcal{X}-AP(w))-P(\bar{\omega})(I_\mathcal{X}-A^2)P(w)]
(I_\mathcal{X}-AP(w))^{-1}B\\
& & =
B^*(I_\mathcal{X}-P(\bar{\omega})A)^{-1}(I_\mathcal{X}-P(\bar{\omega})P(w))
(I_\mathcal{X}-AP(w))^{-1}B\\
& & =\sum\limits_{k=1}^N(1-\overline{\omega_k}w_k)B^*
(I_\mathcal{X}-P(\bar{\omega})A)^{-1}P_{\mathcal{X}_k}(I_\mathcal{X}-AP(w))^{-1}B.
\end{eqnarray*}
Analogously,
$$\mathcal{F}(w)-\mathcal{F}(\omega)^*=
\sum\limits_{k=1}^N(w_k-\overline{\omega_k})B^*
(I_\mathcal{X}-P(\bar{\omega})A)^{-1}P_{\mathcal{X}_k}(I_\mathcal{X}-AP(w))^{-1}B.$$
Thus, from (\ref{eq:F+}) and (\ref{eq:F-}) we obtain that $F$
satisfies (\ref{eq:AH+}) and (\ref{eq:AH-}) with
$\Xi_k(w,\omega)=\xi_k(\omega)^*\xi_k(w)$, where
$$\xi_k(w)=\sqrt{2}P_{\mathcal{X}_k}(I_\mathcal{X}-AP(w))^{-1}B
(I_\mathcal{U}-\mathcal{F}(w))^{-1},\quad w\in\mathbb{D}^N,\
k=1,\ldots,N.$$ Since for $z_k\in\Pi,\zeta_k\in\Pi$ one has
$$1-\frac{\overline{\zeta_k}-1}{\overline{\zeta_k}+1}\cdot
\frac{z_k-1}{z_k+1}=\frac{2(z_k+\overline{\zeta_k})}{(\overline{\zeta_k}+1)(z_k+1)},
\quad
\frac{z_k-1}{z_k+1}-\frac{\overline{\zeta_k}-1}{\overline{\zeta_k}+1}=
\frac{2(z_k-\overline{\zeta_k})}{(\overline{\zeta_k}+1)(z_k+1)},$$
by setting $w_k=\frac{z_k-1}{z_k+1}$ and
$\omega_k=\frac{\zeta_k-1}{\zeta_k+1},\ k=1,\ldots,N$, in
(\ref{eq:AH+}) and (\ref{eq:AH-}), we get for
$f(z)=F(\frac{z_1-1}{z_1+1},\ldots,\frac{z_N-1}{z_N+1})$ the
identities (\ref{eq:b+}) and (\ref{eq:b-}) with
$\Phi_k(z,\zeta)=\varphi_k(\zeta)^*\varphi_k(z),\
(z,\zeta)\in\Pi^N\times\Pi^N$, and
\begin{eqnarray*}
\varphi_k(z)=\frac{\sqrt{2}}{z_k+1}
\xi_k\left(\frac{z_1-1}{z_1+1},\ldots,\frac{z_N-1}{z_N+1}\right)\hspace{6.2cm}\\
=
\frac{2}{z_k+1}P_{\mathcal{X}_k}\left(I_\mathcal{X}-A\sum\limits_{k=1}^N
\frac{z_k-1}{z_k+1}P_{\mathcal{X}_k}\right)^{-1}B\left(I_\mathcal{U}-\mathcal{F}\left
(\frac{z_1-1}{z_1+1},\ldots,\frac{z_N-1}{z_N+1}\right)\right)^{-1},
\end{eqnarray*}
$k=1,\ldots,N$. Thus, finally we get $\mathcal{F}=\mathcal{C}(f)$
where $f\in\mathcal{B}_N(\mathcal{U})$. The proof is complete.
\end{proof}

\section{De-homogenization}\label{sec:dehom}
In this section we establish a one-to-one correspondence between
the classes $\mathcal{P}_N(\mathcal{U})$,
$\mathcal{B}_N(\mathcal{U})$ and certain classes of
non-homogeneous holomorphic functions of $N-1$ variables.
\begin{thm}\label{thm:dehom-p}
Let $f\in\mathcal{P}_N(\mathcal{U})$. Then the function
\begin{equation}\label{eq:g}
    g(z_1,\ldots,z_{N-1}):=f(z_1,\ldots,z_{N-1},1)
\end{equation}
is holomorphic in the domain
$\Omega_{N-1}^+:=\bigcup_{\lambda\in\mathbb{T}\cap\Pi}(\lambda\Pi)^{N-1}
\subset\mathbb{C}^{N-1}$ and has the properties:
\begin{description}
    \item[(i)] $g(\overline{z'})=g(z')^*,\quad z'\in\Omega_{N-1}^+$;
    \item[(ii)] $\mbox{Re}\left[z_Ng\left(\frac{z_1}{z_N},
    \ldots,\frac{z_{N-1}}{z_N}\right)\right]\geq 0,\quad
    z=(z_1,\ldots,z_{N-1},z_N)\in\Pi^N$.
\end{description}
Conversely, let an $L(\mathcal{U})$-valued function
$g(z')=g(z_1,\ldots,z_{N-1})$ be holomorphic in the domain
$\Omega_{N-1}^+$ and satisfy conditions (i) and (ii). Then the
function
\begin{equation}\label{eq:f}
    f(z_1,\ldots,z_{N-1},z_N):=z_Ng\left(\frac{z_1}{z_N},
    \ldots,\frac{z_{N-1}}{z_N}\right)
\end{equation}
is correctly defined on $\Omega_N$ and belongs to the class
$\mathcal{P}_N(\mathcal{U})$.
\end{thm}
\begin{proof}
The function $g$ defined by (\ref{eq:g}) is holomorphic in the
domain
\begin{eqnarray*}
& & \{ z'=(z_1,\ldots,z_{N-1})\in\mathbb{C}^{N-1}:\
(z_1,\ldots,z_{N-1},1)\in\Omega_N\}\\
&=&\{ z'\in\mathbb{C}^{N-1}:\ \exists\lambda\in\mathbb{T}:\
(z_1,\ldots,z_{N-1},1)\in(\lambda\Pi)^N\}\\
&=&\{ z'\in\mathbb{C}^{N-1}:\ \exists\lambda\in\mathbb{T}\cap\Pi:\
z'\in(\lambda\Pi)^{N-1}\} =\Omega_{N-1}^+,
\end{eqnarray*}
since $1\in\lambda\Pi$ means $\lambda\in\Pi$. Clearly, (i) is
valid for $g$ since (\ref{eq:gsym}) is valid for $f$. For any
$z=(z_1,\ldots,z_{N-1},z_N)\in\Pi^N$, due to (\ref{eq:ghom}) and
(\ref{eq:gpos}) for $f$, one has
\begin{eqnarray*}
\mbox{Re}\left[z_Ng\left(\frac{z_1}{z_N},
    \ldots,\frac{z_{N-1}}{z_N}\right)\right] &=&\mbox{Re}\left[z_Nf\left(\frac{z_1}{z_N},
    \ldots,\frac{z_{N-1}}{z_N},1\right)\right]\\
    &=& \mbox{Re}f(z_1,\ldots,z_{N-1},z_N)\geq 0.
\end{eqnarray*}
Thus, (ii) is valid, too.

Conversely, let an $L(\mathcal{U})$-valued function
$g(z')=g(z_1,\ldots,z_{N-1})$ be holomorphic in $\Omega_{N-1}^+$
and satisfy (i) and (ii). Then the function $f$ is correctly
defined by (\ref{eq:f}). Indeed, for any $z\in\Omega_N$ there is a
$\lambda\in\mathbb{T}$ such that $z\in(\lambda\Pi)^N$. Hence,
$\left(\frac{z_1}{z_N},\ldots,\frac{z_{N-1}}{z_N}\right)\in\left(
e^{i(\arg \lambda-\arg
    z_N)}\Pi\right)^{N-1}$. Since $\arg z_N\in (\arg\lambda-\frac{\pi}{2},\arg\lambda
    +\frac{\pi}{2})$, one has $\arg \lambda - \arg z_N\in
    (-\frac{\pi}{2},\frac{\pi}{2})$. Therefore, $\left(\frac{z_1}{z_N},
    \ldots,\frac{z_{N-1}}{z_N}\right)\in\Omega_{N-1}^+$. Moreover,
    $f$ is holomorphic in $\Omega_N$. Properties
    (\ref{eq:ghom})--(\ref{eq:gsym}) of $f$ are easily verified.
    Thus, $f\in\mathcal{P}_N(\mathcal{U})$. The proof is complete.
\end{proof}
\begin{thm}\label{thm:dehom-b}
Let $f\in\mathcal{B}_N(\mathcal{U})$. Then the function $g$
defined by (\ref{eq:g}) is holomorphic in the domain
$\Omega_{N-1}^+$ and has the properties:
\begin{description}
    \item[(i)] $g(\overline{z'})=g(z')^*,\quad z'\in\Omega_{N-1}^+$;
    \item[(ii)] for any $\mathbf{R}\in\mathcal{A}^N$ such that
    $R_k\in L(\mathcal{H}),\
    k=1,\ldots,N$, where $\mathcal{H}$ is some Hilbert space, the operator
     $g(R_N^{-1}R_1,\ldots,R_N^{-1}R_{N-1})$ is correctly defined
     and
     $\mbox{Re}[(I_\mathcal{U}\otimes
     R_N)g(R_N^{-1}R_1,\ldots,R_N^{-1}R_{N-1})]\geq 0.$
\end{description}
Conversely, let an $L(\mathcal{U})$-valued function
$g(z')=g(z_1,\ldots,z_{N-1})$ be holomorphic in the domain
$\Omega_{N-1}^+$ and satisfy conditions (i) and (ii). Then the
function $f$ is correctly defined by (\ref{eq:f}) and belongs to
the class $\mathcal{B}_N(\mathcal{U})$.
\end{thm}
\begin{proof}
Since
$f\in\mathcal{B}_N(\mathcal{U})\subset\mathcal{P}_N(\mathcal{U})$,
by Theorem~\ref{thm:dehom-p} the function $g$  is holomorphic in
$\Omega_{N-1}^+$ and satisfies (i). Since the function
$g\left(\frac{z_1}{z_N},
    \ldots,\frac{z_{N-1}}{z_N}\right)=f\left(\frac{z_1}{z_N},
    \ldots,\frac{z_{N-1}}{z_N},1\right)$
    is holomorphic in $z=(z_1,\ldots,z_N)\in\Omega_N$,
    the operator $g(R_N^{-1}R_1,\ldots,R_N^{-1}R_{N-1})$
is correctly defined for any
$\mathbf{R}=(R_1,\ldots,R_N)\in\mathcal{A}^N$ (see
Section~\ref{sec:b-calc}). Moreover,
 \begin{eqnarray*}
\mbox{Re}[(I_\mathcal{U}\otimes
     R_N)g(R_N^{-1}R_1,\ldots,R_N^{-1}R_{N-1})] &=& \\
 \mbox{Re}[(I_\mathcal{U}\otimes
     R_N)f(R_N^{-1}R_1,\ldots,R_N^{-1}R_{N-1},I_\mathcal{H})]
         &=& \mbox{Re}f(R_1,\ldots,R_{N-1},R_N)\geq 0,
\end{eqnarray*}
due to properties (i) and (ii) of $f$ in Corollary~\ref{cor:bes}.
Thus $g$ has property (ii) of the present theorem.

Conversely, let an $L(\mathcal{U})$-valued function $g(z')$ be
holomorphic in $\Omega_{N-1}^+$ and satisfy conditions (i) and
(ii) of this theorem. In particular, $g$ satisfies conditions (i)
and (ii) of Theorem~\ref{thm:dehom-p}. Then the function $f$ is
correctly defined by (\ref{eq:f}) and belongs to the class
$\mathcal{P}_N(\mathcal{U})$. According to
Corollary~\ref{cor:bes}, the only thing we have to check is that
$\mbox{Re}f(\mathbf{R})\geq 0$ for any
$\mathbf{R}\in\mathcal{A}^N$. The latter follows easily:
$$\mbox{Re}f(R_1,\ldots,R_{N-1},R_N)=\mbox{Re}[(I_\mathcal{U}\otimes
     R_N)g(R_N^{-1}R_1,\ldots,R_N^{-1}R_{N-1})]\geq 0.$$
     The proof is complete.
\end{proof}

\section{The ``real" case}\label{sec:real}
Let us introduce the operator analogues of real matrices, and
operator-valued analogues of real matrix-valued functions. To this
end, first of all we define an anti-linear operator which
generalizes the complex conjugation in $\mathbb{C}^n$.

The operator $\iota$ on a Hilbert space $\mathcal{U}$ is called an
\emph{anti-unitary involution} (AUI) if the following conditions
are satisfied:
\begin{equation}\label{eq:inv}
    \iota^2=I_\mathcal{U};
\end{equation}
\begin{equation}\label{eq:anti}
   \langle\iota u_1,\iota
    u_2\rangle=\langle u_2,u_1\rangle,\quad u_1\in\mathcal{U},\
    u_2\in\mathcal{U}.
\end{equation}
\begin{prop}\label{prop:iota}
An AUI $\iota$ has the following properties:
\begin{description}
    \item[(i)] $\iota$ is \textbf{additive}, i.e., for any $u_1\in\mathcal{U},\
    u_2\in\mathcal{U}$ one has $\iota(u_1+u_2)=\iota u_1+\iota
    u_2$;
    \item[(ii)] $\iota$ is \textbf{anti-homogeneous}, i.e., for any $u\in\mathcal{U},\
    \alpha\in\mathbb{C}$ one has $\iota(\alpha u)=\bar{\alpha}\iota
    u$;
    \item[(iii)] the operators $\pi_\pm
    :=\frac{I_\mathcal{U}\pm\iota}{2}$ are idempotents, their
    ranges $\pi_\pm\mathcal{U}$ are closed in $\mathcal{U}$, and
    $\pi_+\mathcal{U}\dotplus\pi_-\mathcal{U}=\mathcal{U}$.
\end{description}
\end{prop}
\begin{proof}
Let $u_1\in\mathcal{U},\ u_2\in\mathcal{U}$. Then
\begin{eqnarray*}
\lefteqn{\|\iota(u_1+u_2)-\iota u_1-\iota
    u_2\|^2=\|\iota(u_1+u_2)\|^2+\|\iota u_1\|^2+\|\iota
    u_2\|^2}\\
&-& 2\mbox{Re}\langle\iota(u_1+u_2),\iota u_1+\iota
    u_2\rangle +2\mbox{Re}\langle\iota u_1,\iota u_2\rangle \\
&= &\| u_1+u_2\|^2+\| u_1\|^2+\| u_2\|^2-2\mbox{Re}\langle
    u_1,u_1+u_2\rangle -2\mbox{Re}\langle u_2,u_1+u_2\rangle \\
&+& 2\mbox{Re}\langle u_2,u_1\rangle = \| u_1\|^2+\| u_2\|^2-\|
u_1+u_2\|^2+2\mbox{Re}\langle
    u_2,u_1\rangle =0.
\end{eqnarray*}
Therefore, $\iota(u_1+u_2)-\iota u_1-\iota u_2=0$, that proves
(i). Next, for any $u\in\mathcal{U},\ \alpha\in\mathbb{C}$ one has
\begin{eqnarray*}
\lefteqn{\|\iota(\alpha u)-\bar{\alpha}\iota u\|^2=\|\iota(\alpha
u)\|^2-2\mbox{Re}\langle\iota(\alpha u),\bar{\alpha}\iota u\rangle
+|\alpha|^2\|\iota u\|^2}\\
&=&\|\alpha u\|^2-2\mbox{Re}(\alpha\langle\iota(\alpha u),\iota
u\rangle )
+|\alpha|^2\| u\|^2\\
&=& 2|\alpha|^2\| u\|^2-2\mbox{Re}(\alpha\langle u,\alpha u\rangle
)= 2|\alpha|^2\| u\|^2- 2|\alpha|^2\| u\|^2=0,
\end{eqnarray*}
that proves (ii). Since
$$\pi_\pm^2=\frac{(I_\mathcal{U}\pm\iota)^2}{4}=\frac{(I_\mathcal{U}\pm 2\iota
+\iota^2)}{4}=\frac{I_\mathcal{U}\pm\iota}{2}=\pi_\pm ,$$ the
operators $\pi_\pm$ are idempotents. Since $\iota$ preserves the
norm, $\iota$ is continuous, thus $\pi_\pm$ are also continuous
operators. The latter means that if $\pi_\pm u_j\rightarrow y$ as
$j\rightarrow\infty$ then $\pi_\pm u_j=\pi_\pm(\pi_\pm
u_j)\rightarrow\pi_\pm y$, i.e., $y=\pi_\pm
y\in\pi_\pm\mathcal{U}$. Therefore, $\pi_\pm\mathcal{U}$ are
closed lineals in $\mathcal{U}$.

If $\pi_+u_1=\pi_-u_2$ then $\pi_+^2u_1=\pi_+\pi_-u_2$. Since
$\pi_+^2=\pi_+$ and
$\pi_+\pi_-=\frac{(I_\mathcal{U}+\iota)(I_\mathcal{U}-\iota)}{4}=
\frac{I_\mathcal{U}-\iota^2}{4}=0$, one has $\pi_+u_1=\pi_-u_2=0$,
i.e., $\pi_+\mathcal{U}\cap\pi_-\mathcal{U}=\{ 0\}$. Since for any
$u\in\mathcal{U}$ one has $u=\pi_+u+\pi_-u$, and $\pi_\pm
u\in\pi_\pm\mathcal{U}$, we get
$\pi_+\mathcal{U}\dotplus\pi_-\mathcal{U}=\mathcal{U}$. The proof
is complete.
\end{proof}
\begin{ex}\label{ex:conj}
Let $\mathcal{U}=\mathbb{C}^n$. For
$u=\mbox{col}(u_1,\ldots,u_n)\in\mathbb{C}^n$ denote
$\bar{u}:=\mbox{col}(\overline{u_1},\ldots,\overline{u_n})$. Set
$\iota u:=\bar{u}$, i.e, $\iota$ is the complex conjugation in
$\mathbb{C}^n$. Then, clearly, $\iota$ is an AUI and
$$\pi_+u=\mbox{Re}\,u=\mbox{col}(\mbox{Re}\,u_1,\ldots,\mbox{Re}\,u_n),\
\pi_-u=i\mbox{Im}\,u=i\mbox{col}(\mbox{Im}\,u_1,\ldots,\mbox{Im}\,u_n).$$
\end{ex}
The operator $A\in L(\mathcal{U})$ is called \emph{$\iota$-real}
(resp., \emph{$\iota$-symmetric}) if $\iota A=A\iota$ (resp.,
$\iota A=A^*\iota$).
\begin{ex}\label{ex:real}
Let  $\mathcal{U}=\mathbb{C}^n$ and $\iota u=\bar{u}$, as in
Example~\ref{ex:conj}. Then the operator $A\in L(\mathcal{U})$ is
$\iota$-real (resp., $\iota$-symmetric) if and only if its matrix
in the standard basis $\{ e_k\}_{k=1}^n$ has real entries (resp.,
is symmetric, i.e., $A^T=A$). Indeed, in the first case
\begin{eqnarray*}
a_{kj} &=& \langle Ae_j,e_k\rangle =\langle\iota e_k,\iota
Ae_j\rangle =\langle\iota e_k,A\iota e_j\rangle\\
&=& \langle e_k,Ae_j\rangle =\overline{\langle Ae_j,e_k\rangle
}=\overline{a_{kj}},\quad k,j=1,\ldots,n;
\end{eqnarray*}
in the second case
\begin{eqnarray*}
a_{kj} &=& \langle Ae_j,e_k\rangle =\langle\iota e_k,\iota
Ae_j\rangle =\langle\iota e_k,A^*\iota e_j\rangle\\
&=& \langle e_k,A^*e_j\rangle =\langle Ae_k,e_j\rangle
=a_{jk},\quad k,j=1,\ldots,n.
\end{eqnarray*}
\end{ex}
\begin{lem}\label{lem:iota}
Let $\iota=\iota_\mathcal{U}$ be an AUI on a Hilbert space
$\mathcal{U}$, and $A\in L(\mathcal{U})$. For the following three
conditions, any two of them imply the third one:
\begin{description}
    \item[(i)] $\iota A=A\iota$;
    \item[(ii)] $\iota A=A^*\iota$;
    \item[(iii)] $A=A^*$.
\end{description}
\end{lem}
\begin{proof}
(i)\&(ii)$\Rightarrow$(iii). If $\iota A=A\iota$ and $\iota
A=A^*\iota$ then $A\iota =A^*\iota$. Therefore,
$A\iota^2=A^*\iota^2$, i.e., $A=A^*$.

Implications (i)\&(iii)$\Rightarrow$(ii) and
(ii)\&(iii)$\Rightarrow$(i) are obvious.
\end{proof}
Let $\iota$ be an AUI on a Hilbert space $\mathcal{U}$, and
$\Omega\subset\mathbb{C}^N$ be a domain invariant under the
complex conjugation: $\overline{\Omega}=\Omega$. For a function
$f:\Omega\rightarrow L(\mathcal{U})$ set $f^\sharp(z)=\iota
f(\bar{z})\iota,\ z\in\Omega$. A function $f:\Omega\rightarrow
L(\mathcal{U})$ is called \emph{$\iota$-real} if for any
$z\in\Omega$ one has $f^\sharp(z)=f(z)$.
\begin{ex}\label{ex:funct}
Let  $\mathcal{U}=\mathbb{C}^n$ and $\iota u=\bar{u}$, as in
Examples~\ref{ex:conj} and \ref{ex:real}, and
$\Omega\subset\mathbb{C}^N$ be a domain satisfying
$\overline{\Omega}=\Omega$. Then $\iota$-real
$L(\mathcal{U})$-valued functions are those matrix functions whose
values in the standard basis $\{ e_k\}_{k=1}^n$ satisfy the
condition $f(\bar{z})=\overline{f(z)}$. Indeed,
\begin{eqnarray*}
[f(\bar{z})]_{kj} &=& \langle f(\bar{z})e_j,e_k\rangle =\langle
f(\bar{z})\iota e_j,\iota
e_k\rangle =\langle\iota^2 e_k,\iota f(\bar{z})\iota e_j\rangle\\
&=& \langle e_k,f(z)e_j\rangle =\overline{\langle
f(z)e_j,e_k\rangle }=\overline{[f(z)]_{kj}},\quad k,j=1,\ldots,n.
\end{eqnarray*}
\end{ex}
Let $\mathcal{U}$ be a Hilbert space and $\iota=\iota_\mathcal{U}$
be an AUI on $\mathcal{U}$. Denote by
$\iota\mathbb{R}\mathcal{P}_N(\mathcal{U})$ and
$\iota\mathbb{R}\mathcal{B}_N(\mathcal{U})$ the subclasses of
$\mathcal{P}_N(\mathcal{U})$ and  $\mathcal{B}_N(\mathcal{U})$,
respectively, consisting of $\iota$-real functions. Clearly,
$\iota\mathbb{R}\mathcal{B}_N(\mathcal{U})\subset
\iota\mathbb{R}\mathcal{P}_N(\mathcal{U})$.
\begin{prop}\label{prop:Taylor}
The Taylor coefficients of functions from
$\iota\mathbb{R}\mathcal{P}_N(\mathcal{U})$, and therefore,
functions from $\iota\mathbb{R}\mathcal{B}_N(\mathcal{U})$ satisfy
conditions (i)--(iii) of Lemma~\ref{lem:iota}.
\end{prop}
\begin{proof}
According to Lemma \ref{lem:iota}, it is sufficient to verify any
two of its conditions, e.g., (i) and (iii), for the Taylor
coefficients of an arbitrary function
$f\in\iota\mathbb{R}\mathcal{P}_N(\mathcal{U})$. It follows from
Proposition~\ref{prop:iota} that
$$f^\sharp(z)=\iota f(\bar{z})\iota=\sum_{t\in\mathbb{Z}^N_+}\iota\widehat{f_t}
(\bar{z}-\overline{z^\circ})^t\iota
=\sum_{t\in\mathbb{Z}^N_+}\iota\widehat{f_t}\iota(z-z^\circ)^t$$
in some neighborhood $\Gamma$ of $z^\circ\in\Omega_N$. Since $f$
is $\iota$-real,
$$f^\sharp(z)=f(z)=\sum_{t\in\mathbb{Z}^N_+}\widehat{f_t}(z-z^\circ)^t$$
in some neighborhood $\Gamma'$ of $z^\circ\in\Omega_N$. By the
uniqueness of Taylor's expansion in the neighborhood
$\Gamma\cap\Gamma'$ of $z^\circ$, we get
$\iota\widehat{f_t}\iota=\widehat{f_t}$ and
$\iota\widehat{f_t}=\widehat{f_t}\iota$, for any
$t\in\mathbb{Z}^N_+$, i.e., condition (i) is satisfied  for all
the Taylor coefficients of $f$. Since
$$f(\bar{z})^*=\left(\sum_{t\in\mathbb{Z}^N_+}\widehat{f_t}
(\bar{z}-\overline{z^\circ})^t\right)^*
=\sum_{t\in\mathbb{Z}^N_+}\widehat{f_t}^*(z-z^\circ)^t$$ in
$\Gamma$, and
$$f(\bar{z})^*=f(z)=\sum_{t\in\mathbb{Z}^N_+}\widehat{f_t}(z-z^\circ)^t$$
in $\Gamma'$, by the uniqueness of Taylor's expansion in the
neighborhood $\Gamma\cap\Gamma'$ of $z^\circ$, we get
$\widehat{f_t}^*=\widehat{f_t}$ for any $t\in\mathbb{Z}^N_+$, i.e,
condition (iii) is also satisfied  for all the Taylor coefficients
of $f$.
\end{proof}
\begin{thm}\label{thm:real}
Let $f$ be a holomorphic $L(\mathcal{U})$-valued function on
$\Omega_N$, and $\iota=\iota_\mathcal{U}$ be an AUI on a Hilbert
space $\mathcal{U}$. The following conditions are equivalent:
\begin{description}
    \item[(i)] $f\in\iota\mathbb{R}\mathcal{B}_N(\mathcal{U})$;
    \item[(ii)] there exist a representation  (\ref{eq:gsc}) of
    $f$ and an AUI $\iota_\mathcal{H}$ on $\mathcal{H}$ such that
    the function $A(z)$ given by  (\ref{eq:glp}) is
    $(\iota_\mathcal{U}\oplus\iota_\mathcal{H})$-real;
    \item[(iii)] there exists a representation  (\ref{eq:kerid})
    of $f$ for which the holomorphic PSD kernels $\Phi_k(z,\zeta),\ k=1,\ldots,N$, on
    $\Pi^N\times\Pi^N$ are $\iota_\mathcal{U}$-real;
    \item[(iv)] there exist an Agler unitary colligation
    $\alpha=(N;U;\mathcal{X}=\bigoplus_{k=1}^N\mathcal{X}_k,
    \mathcal{U},\mathcal{U})$ with $U=U^*$, the corresponding
    representation (\ref{eq:Ag}) of $\mathcal{F}=\mathcal{C}(f)$
    (the latter is given by (\ref{eq:Cayley-val})) and an AUI
    $\iota_\mathcal{X}$ on $\mathcal{X}$ which commutes with the
    orthogonal projectors $P_k=P_{\mathcal{X}_k},\ k=1,\ldots,N$,
    such that the operator $U$ is
    $(\iota_\mathcal{X}\oplus\iota_\mathcal{U})$-real.
\end{description}
\end{thm}
\begin{proof}
(i)$\Rightarrow$(iii). By Theorem~\ref{thm:ker} there exists a
representation (\ref{eq:kerid}) of $f$ with holomorphic PSD
kernels $\Phi_k(z,\zeta),\ k=1,\ldots,N$, on
    $\Pi^N\times\Pi^N$. Define
$$\Phi^\circ_k(z,\zeta):=\frac{\Phi_k(z,\zeta)+\Phi_k^\sharp(z,\zeta)}{2},
\quad (z,\zeta)\in\Pi^N\times\Pi^N,\ k=1,\ldots,N.$$ For any
$u_1\in\mathcal{U},\ u_2\in\mathcal{U},\ k=1,\ldots,N$, one has
\begin{eqnarray*}
\langle\Phi_k^\sharp(z,\overline{\zeta})u_1,u_2\rangle &=&
\langle\iota\Phi_k(\bar{z},\zeta)\iota u_1,u_2\rangle=
\langle\iota u_2,\Phi_k(\bar{z},\zeta)\iota u_1\rangle
\\
&=& \langle\Phi_k(\bar{z},\zeta)^*\iota u_2,\iota u_1\rangle
=\langle\Phi_k(\zeta,\bar{z})\iota u_2,\iota u_1\rangle,
\end{eqnarray*}
which is, clearly, a holomorphic function in
$(z,\zeta)\in\Pi^N\times\Pi^N$ (we used the property (\ref{eq:K1})
of PSD kernels). Thus, $\Phi_k^\sharp(z,\overline{\zeta})$ and,
therefore, $\Phi_k^\circ(z,\overline{\zeta})$ are holomorphic
operator-valued functions on $\Pi^N\times\Pi^N$. Since for any
$u_1\in\mathcal{U},\ u_2\in\mathcal{U}$
\begin{eqnarray*}
\langle\Phi_k^\sharp(z,\zeta)u_1,u_2\rangle &=&
\langle\Phi_k(\overline{\zeta},\bar{z})\iota u_2,\iota u_1\rangle
=\langle u_1,\iota\Phi_k(\bar{\zeta},\bar{z})\iota u_2\rangle \\
&=& \langle u_1,\Phi_k^\sharp(\zeta,z) u_2\rangle,
\end{eqnarray*}
one has $\Phi_k^\sharp(\zeta,z)=\Phi_k^\sharp(z,\zeta)^*$ and,
therefore, $\Phi_k^\circ(\zeta,z)=\Phi_k^\circ(z,\zeta)^*$ for all
$(z,\zeta)\in\Pi^N\times\Pi^N,\ k=1,\ldots,N$, i.e.,
$\Phi_k^\circ(z,\zeta)$ satisfies condition (\ref{eq:K1}).

 Let
$m\in\mathbb{N},\ \{ z^{(\mu)}\}_{\mu=1}^m\subset\Pi^N,\ \{
u^{(\mu)}\}_{\mu=1}^m\subset\mathcal{U}$. Then for $k=1,\ldots,N$
one has
\begin{eqnarray*}
\sum_{\mu=1}^m\sum_{\nu=1}^m\langle\Phi_k^\sharp(z^{(\mu)},z^{(\nu)})u^{(\mu)},
u^{(\nu)}\rangle &=&
\sum_{\mu=1}^m\sum_{\nu=1}^m\langle\iota\Phi_k(\overline{z^{(\mu)}},
\overline{z^{(\nu)}})\iota u^{(\mu)}, u^{(\nu)}\rangle \\
\sum_{\mu=1}^m\sum_{\nu=1}^m\langle\iota
u^{(\nu)},\Phi_k(\overline{z^{(\mu)}}, \overline{z^{(\nu)}})\iota
u^{(\mu)}\rangle &=& \sum_{\mu=1}^m\sum_{\nu=1}^m\langle\iota
u^{(\nu)},\Phi_k(\overline{z^{(\nu)}},
\overline{z^{(\mu)}})^*\iota u^{(\mu)}\rangle \\
&=&
\sum_{\mu=1}^m\sum_{\nu=1}^m\langle\Phi_k(\overline{z^{(\nu)}},
\overline{z^{(\mu)}})\iota u^{(\nu)},\iota u^{(\mu)}\rangle\geq 0,
\end{eqnarray*}
by virtue of (\ref{eq:K2}) for $\Phi_k(z,\zeta)$. Thus,
$\Phi_k^\sharp(z,\zeta)$ and, therefore, $\Phi_k^\circ(z,\zeta)$
satisfy condition (\ref{eq:K2}). Finally, we have proved that
$\Phi_k^\sharp(z,\zeta)$ and, therefore, $\Phi_k^\circ(z,\zeta)$
are holomorphic PSD kernels on $\Pi^N\times\Pi^N$.

Since $f^\sharp(z)=\iota f(\bar{z})\iota=f(z)$, the identity
(\ref{eq:kerid}) together with properties (i) and (ii) in
Proposition~\ref{prop:iota} of $\iota$ imply
\begin{eqnarray*}
f(z) &=&
\iota\left(\sum_{k=1}^N\overline{z_k}\Phi_k(\bar{z},\bar{\zeta})\right)\iota=
\sum_{k=1}^Nz_k\iota\Phi_k(\bar{z},\bar{\zeta})\iota \\
&=& \sum_{k=1}^Nz_k\Phi_k^\sharp(z,\zeta),\quad
(z,\zeta)\in\Pi^N\times\Pi^N.
\end{eqnarray*}
Therefore, $f(z)=\sum_{k=1}^Nz_k\Phi_k^\circ(z,\zeta),\quad
(z,\zeta)\in\Pi^N\times\Pi^N.$ Moreover,
\begin{eqnarray*}
\Phi_k^{\circ\sharp}(z,\zeta)=\frac{\Phi_k^\sharp(z,\zeta)+\Phi_k^{\sharp\sharp}
(z,\zeta)}{2}=\frac{\Phi_k^\sharp(z,\zeta)+\Phi_k
(z,\zeta)}{2}=\Phi_k^\circ(z,\zeta).
\end{eqnarray*}
Thus, (\ref{eq:kerid}) holds with $\iota_\mathcal{U}$-real
holomorphic PSD kernels $\Phi_k^\circ(z,\zeta),\ k=1,\ldots,N$,
and (iii) follows from (i).

(iii)$\Rightarrow$(ii). Recall that holomorphic PSD kernels
$\Phi_k(z,\zeta)$ admit factorizations
$\Phi_k(z,\zeta)=\varphi_k(\zeta)^*\varphi_k(z)$, where
$\varphi_k$ are holomorphic $L(\mathcal{U,M}_k)$-valued functions
on $\Pi^N$, and $\mathcal{M}_k$ are auxiliary Hilbert spaces,
$k=1,\ldots,N$. These $\mathcal{M}_k$'s and $\varphi_k$'s can be
determined in the following way. Set
$\widetilde{\Phi_k}(z,\zeta):=\Phi_k(\bar{z},\bar{\zeta})^*$. It
is easy to convince oneself that $\widetilde{\Phi_k}(z,\zeta)$ are
holomorphic PSD kernels on $\Pi^N\times\Pi^N$. Then set
$\mathcal{M}_k:=\mathcal{H}_{\widetilde{\Phi_k}}$, where
$\mathcal{H}_{\widetilde{\Phi_k}}$ denotes the Hilbert space with
the reproducing kernel $\widetilde{\Phi_k}(z,\zeta)$, which is
obtained by completion of the linear span of functions of the form
$\widetilde{\Phi_k}(\cdot,\zeta)u,\ \zeta\in\Pi^N,\
u\in\mathcal{U}$, with respect to the inner product
$$\langle\widetilde{\Phi_k}(\cdot,\zeta)u,\widetilde{\Phi_k}(\cdot,\zeta')u'
\rangle:= \langle\widetilde{\Phi_k}(\zeta',\zeta)u,u'\rangle.$$
Clearly, $\widetilde{\Phi_k}(z,\zeta)=H_k(z)H_k(\zeta)^*,\
(z,\zeta)\in\Pi^N\times\Pi^N$, where $H_k$ is a holomorphic
$L(\mathcal{M}_k,\mathcal{U})$-valued function on $\Pi^N$ defined
by $H_k(z)^*u:=\widetilde{\Phi_k}(\cdot,z)u,\ u\in\mathcal{U}$.
Set $\varphi_k(z):=\widetilde{H_k}(z)=H_k(\overline{z})^*$. Then
for $u\in\mathcal{U}$ one has
$\varphi_k(z)u=\widetilde{\Phi_k}(\cdot,\bar{z})u$, and
$$\Phi_k(z,\zeta)=\widetilde{\Phi_k}(\bar{z},\bar{\zeta})^*=H_k(\bar{\zeta})
H_k(\bar{z})^*=\varphi_k(\zeta)^*\varphi_k(z).$$ Define the AUI
$\iota_{\mathcal{M}_k}:\mathcal{M}_k\rightarrow\mathcal{M}_k$ on
generating vectors by
\begin{eqnarray}
(\iota_{\mathcal{M}_k}[\widetilde{\Phi_k}(\cdot,\zeta)u])(z) &:=&
\iota_\mathcal{U}\widetilde{\Phi_k}(\bar{z},\zeta)u
\nonumber \\
&(=&
\iota_\mathcal{U}\Phi_k(z,\bar{\zeta})^*u=\iota_\mathcal{U}\Phi_k(\bar{\zeta},z)u
=\Phi_k(\zeta,\bar{z})
\iota_\mathcal{U}u \nonumber \\
&=& \widetilde{\Phi_k}(\bar{\zeta},z)^*\iota_\mathcal{U}u=
\widetilde{\Phi_k}(z,\bar{\zeta})\iota_\mathcal{U}u),\
(z,\zeta)\in\Pi^N\times\Pi^N,u\in\mathcal{U}.\label{eq:iota-m}
\end{eqnarray}
This definition is correct. Indeed,
\begin{eqnarray*}
\langle\iota_{\mathcal{M}_k}\widetilde{\Phi_k}(\cdot,\zeta)u,\iota_{\mathcal{M}_k}
\widetilde{\Phi_k}(\cdot,\zeta')u' \rangle =
\langle\widetilde{\Phi_k}(\cdot,\bar{\zeta})\iota_\mathcal{U}u,
\widetilde{\Phi_k}(\cdot,\overline{\zeta'})\iota_\mathcal{U}u'
\rangle
=\langle\widetilde{\Phi_k}(\overline{\zeta'},\overline{\zeta})\iota_\mathcal{U}u,
\iota_\mathcal{U}
u'\rangle \\
= \langle
u',\iota_\mathcal{U}\widetilde{\Phi_k}(\overline{\zeta'},\overline{\zeta})
\iota_\mathcal{U}u\rangle = \langle
u',\widetilde{\Phi_k}(\zeta',\zeta)u\rangle = \langle
u',\widetilde{\Phi_k}(\zeta,\zeta')^*u\rangle
=\langle \widetilde{\Phi_k}(\zeta,\zeta')u',u\rangle\\
 =\langle\widetilde{\Phi_k}(\cdot,\zeta')u',\widetilde{\Phi_k}(\cdot,\zeta)u
\rangle.
\end{eqnarray*}
Therefore, $\iota_{\mathcal{M}_k}$ preserves the norm of any
vector of the form
$\sum_{\mu=1}^m\widetilde{\Phi_k}(\cdot,\zeta^{(\mu)})u^{(\mu)}$.
The density of such vectors in $\mathcal{M}_k$ implies
$$\langle\iota_{\mathcal{M}_k}m,\iota_{\mathcal{M}_k}m'\rangle=
\langle m',m\rangle\quad m\in \mathcal{M}_k,m'\in\mathcal{M}_k,$$
that is an analogue of (\ref{eq:anti}). Next,
\begin{eqnarray*}
(\iota_{\mathcal{M}_k}^2[\widetilde{\Phi_k}(\cdot,\zeta)u])(z) &=&
(\iota_{\mathcal{M}_k}[\widetilde{\Phi_k}(\cdot,\bar{\zeta})\iota_\mathcal{U}u])(z)
=\widetilde{\Phi_k}(z,\zeta)\iota_\mathcal{U}^2u \\
&=& \widetilde{\Phi_k}(z,\zeta)u, \quad
(z,\zeta)\in\Pi^N\times\Pi^N,u\in\mathcal{U}.
\end{eqnarray*}
Therefore, by the additivity of $\iota_{\mathcal{M}_k}$ (and
$\iota_{\mathcal{M}_k}^2$) and the linearity and continuity
($\iota_{\mathcal{M}_k}$ is norm-preserving!) argument, we get
$\iota_{\mathcal{M}_k}^2=I_{\mathcal{M}_k}$, i.e., an analogue of
(\ref{eq:inv}).

The identities in (\ref{eq:iota-m}) imply
\begin{equation}\label{eq:iota-m'}
    \iota_{\mathcal{M}_k}\varphi_k(z)=\varphi_k(\bar{z})\iota_\mathcal{U},\quad
    z\in\Pi^N.
\end{equation}
Following the sufficiency part of the proof of
Theorem~\ref{thm:ker}, we obtain
$\varphi(z)=\mbox{col}[\begin{array}{ccc}
  \varphi_1(z) & \ldots & \varphi_N(z)
\end{array}]$,
$$\mathcal{H}=\mbox{clos span}_{z\in\Pi^N}\{
[\varphi(z)-\varphi(e)]\mathcal{U}\}\subset\mathcal{M}=
\bigoplus_{k=1}^N\mathcal{M}_k,$$
$$A_k=\left[\begin{array}{cc}
  \varphi(e)^* & 0 \\
  0 & I_\mathcal{H}
\end{array}\right]\kappa^*P_k\kappa\left[\begin{array}{cc}
  \varphi(e) & 0 \\
  0 & I_\mathcal{H}
\end{array}\right]\in L(\mathcal{U}\oplus\mathcal{H})$$
(here $\kappa$ is defined by (\ref{eq:kappa})) such that
(\ref{eq:gsc}) holds for $f$, where $$A_k=\left[\begin{array}{cc}
  a_k & b_k \\
  c_k & d_k
\end{array}\right]\in L(\mathcal{U}\oplus\mathcal{H}),\quad k=1,\ldots,N.$$
For $\psi(z)=\left[\begin{array}{c}
  I_\mathcal{U} \\
  \varphi(z)-\varphi(e)
\end{array}\right]$ one has $\psi(e)\mathcal{U}=\mathcal{U}\oplus\{
0\}$, therefore the linear span of vectors of the form $\psi(z)u,\
z\in\Pi^N,\ u\in\mathcal{U}$, is dense in
$\mathcal{U}\oplus\mathcal{H}$. Set
$\iota_\mathcal{M}:=\bigoplus_{k=1}^N\iota_{\mathcal{M}_k}$. By
virtue of (\ref{eq:iota-m'}), one has
$\iota_\mathcal{M}(\varphi(z)-\varphi(e))u=
(\varphi(\bar{z})-\varphi(e))\iota_\mathcal{U}u\in\mathcal{H},$
hence $\iota_\mathcal{M}\mathcal{H}\subset\mathcal{H}$. Set
$\iota_\mathcal{H}:=\iota_\mathcal{M}|\mathcal{H}$. Clearly,
$\iota_\mathcal{H}$ is an AUI on $\mathcal{H}$, and
$(\iota_\mathcal{U}\oplus\iota_\mathcal{H})\psi(z)=\psi(\bar{z})\iota_\mathcal{U},\
z\in\Pi^N$.

Let us verify
$(\iota_\mathcal{U}\oplus\iota_\mathcal{H})A_k(\iota_\mathcal{U}\oplus
\iota_\mathcal{H}) =A_k$, i.e., the
$(\iota_\mathcal{U}\oplus\iota_\mathcal{H})$-realness of $A_k, \
k=1,\ldots,N$. For any $z\in\Pi^N,\ \zeta\in\Pi^N,\
u\in\mathcal{U},\ u'\in\mathcal{U}$ one has
\begin{eqnarray*}
\left\langle A_k\left[\begin{array}{c}
  I_\mathcal{U} \\
  \varphi(z)-\varphi(e)
\end{array}\right]u,\left[\begin{array}{c}
  I_\mathcal{U} \\
  \varphi(\zeta)-\varphi(e)
\end{array}\right]u'\right\rangle \hspace{5.5cm}\\
= \left\langle\left[\begin{array}{cc}
  \varphi(e)^* & 0 \\
  0 & I_\mathcal{H}
\end{array}\right]\kappa^*P_k\kappa\left[\begin{array}{cc}
  \varphi(e) & 0 \\
  0 & I_\mathcal{H}
\end{array}\right]\left[\begin{array}{c}
  I_\mathcal{U} \\
  \varphi(z)-\varphi(e)
\end{array}\right]u,\left[\begin{array}{c}
  I_\mathcal{U} \\
  \varphi(\zeta)-\varphi(e)
\end{array}\right]u'\right\rangle \\
 = \left\langle P_k\kappa\left[\begin{array}{cc}
  \varphi(e) & 0 \\
  0 & I_\mathcal{H}
\end{array}\right]\left[\begin{array}{c}
  I_\mathcal{U} \\
  \varphi(z)-\varphi(e)
\end{array}\right]u,\kappa\left[\begin{array}{cc}
  \varphi(e) & 0 \\
  0 & I_\mathcal{H}
\end{array}\right]\left[\begin{array}{c}
  I_\mathcal{U} \\
  \varphi(\zeta)-\varphi(e)
\end{array}\right]u'\right\rangle \\
 = \left\langle P_k\kappa\left[\begin{array}{c}
\varphi(e) \\
  \varphi(z)-\varphi(e)
\end{array}\right]u,\kappa\left[\begin{array}{c}
  \varphi(e) \\
  \varphi(\zeta)-\varphi(e)
\end{array}\right]u'\right\rangle =\langle P_k\varphi(z)u,\varphi(\zeta)u'\rangle \\
 = \langle P_k\varphi(z)u,P_k\varphi(\zeta)u'\rangle =
\langle\varphi_k(z)u,\varphi_k(\zeta)u'\rangle;
\end{eqnarray*}
\begin{eqnarray*}
\left\langle
(\iota_\mathcal{U}\oplus\iota_\mathcal{H})A_k(\iota_\mathcal{U}\oplus
\iota_\mathcal{H})\left[\begin{array}{c}
  I_\mathcal{U} \\
  \varphi(z)-\varphi(e)
\end{array}\right]u,\left[\begin{array}{c}
  I_\mathcal{U} \\
  \varphi(\zeta)-\varphi(e)
\end{array}\right]u'\right\rangle \hspace{2.5cm}\\
=
\left\langle(\iota_\mathcal{U}\oplus\iota_\mathcal{H})
\left[\begin{array}{c}
  I_\mathcal{U} \\
  \varphi(\zeta)-\varphi(e)
\end{array}\right]u',A_k(\iota_\mathcal{U}\oplus
\iota_\mathcal{H})\left[\begin{array}{c}
  I_\mathcal{U} \\
  \varphi(z)-\varphi(e)
\end{array}\right]u\right\rangle \\
 = \left\langle \left[\begin{array}{c}
  \iota_\mathcal{U} \\
  \iota_\mathcal{H}(\varphi(\zeta)-\varphi(e))
\end{array}\right]u',A_k\left[\begin{array}{c}
  \iota_\mathcal{U} \\
  \iota_\mathcal{H}(\varphi(z)-\varphi(e))
\end{array}\right]u\right\rangle \\
 =\left\langle \left[\begin{array}{c}
I_\mathcal{U} \\
  \varphi(\overline{\zeta})-\varphi(e)
\end{array}\right]\iota_\mathcal{U}u',A_k\left[\begin{array}{c}
I_\mathcal{U} \\
  \varphi(\overline{z})-\varphi(e)
\end{array}\right]\iota_\mathcal{U}u\right\rangle \\
 =\left\langle A_k\left[\begin{array}{c}
I_\mathcal{U} \\
  \varphi(\overline{\zeta})-\varphi(e)
\end{array}\right]\iota_\mathcal{U}u',\left[\begin{array}{c}
I_\mathcal{U} \\
  \varphi(\overline{z})-\varphi(e)
\end{array}\right]\iota_\mathcal{U}u\right\rangle =
\langle\varphi_k(\overline{\zeta})\iota_\mathcal{U}u',
\varphi_k(\overline{z})\iota_\mathcal{U}u\rangle \\
=\langle\iota_{\mathcal{M}_k}\varphi_k(\zeta)u',\iota_{\mathcal{M}_k}
\varphi_k(z)u\rangle =
\langle\varphi_k(z)u,\varphi_k(\zeta)u'\rangle
\end{eqnarray*}
(in the second chain of calculations we used the result of the
first one and (\ref{eq:iota-m'})). As mentioned above, the linear
span of vectors of the form $\left[\begin{array}{c}
  I_\mathcal{U} \\
  \varphi(z)-\varphi(e)
\end{array}\right]u,\ z\in\Pi^N,\ u\in\mathcal{U}$,
is dense in $\mathcal{U}\oplus\mathcal{H}$. Operators $A_k$ and
$(\iota_\mathcal{U}\oplus\iota_\mathcal{H})A_k(\iota_\mathcal{U}\oplus
\iota_\mathcal{H})$ are continuous and linear (the second operator
is additive because $A_k$ and
$\iota_\mathcal{U}\oplus\iota_\mathcal{H}$ are additive, and
homogeneous because $A_k$ is homogeneous, and
$\iota_\mathcal{U}\oplus\iota_\mathcal{H}$ is anti-homogeneous and
appears twice). Therefore, by comparison the results of the two
chains of calculations above, and linearity and continuity
argument we obtain
$(\iota_\mathcal{U}\oplus\iota_\mathcal{H})A_k(\iota_\mathcal{U}\oplus
\iota_\mathcal{H}) =A_k,\ k=1,\ldots,N$. Thus, (ii) follows from
(iii).

(ii)$\Rightarrow$(i). Let $f$ satisfies (ii). Then
$f\in\mathcal{B}_N(\mathcal{U})$, and the operator-valued linear
function $A(z)$ is
$\iota_\mathcal{U}\oplus\iota_\mathcal{H}$-real, i.e.,
$$\left[\begin{array}{cc}
  \iota_\mathcal{U} & 0 \\
  0 & \iota_\mathcal{H}
\end{array}\right]\left[\begin{array}{cc}
  a(\bar{z}) & b(\bar{z}) \\
  c(\bar{z}) & d(\bar{z})
\end{array}\right]\left[\begin{array}{cc}
  \iota_\mathcal{U} & 0 \\
  0 & \iota_\mathcal{H}
\end{array}\right]=  \left[\begin{array}{cc}
  a(z) & b(z) \\
  c(z) & d(z)
\end{array}\right],\quad z\in\Omega_N.$$
The latter is equivalent to the identities
$$\begin{array}{cc}
  \iota_\mathcal{U}a(\bar{z})\iota_\mathcal{U}=a(z),  &
  \iota_\mathcal{U}b(\bar{z})\iota_\mathcal{H}=b(z), \\
  \iota_\mathcal{H}c(\bar{z})\iota_\mathcal{U}=c(z) &
  \iota_\mathcal{H}d(\bar{z})\iota_\mathcal{H}=d(z),
\end{array}\quad z\in\Omega_N.$$
Since $\iota_\mathcal{H}^2=I_\mathcal{H}$, and
$$ (\iota_\mathcal{H}d(\bar{z})^{-1}\iota_\mathcal{H})\cdot
(\iota_\mathcal{H}d(\bar{z})\iota_\mathcal{H})=
(\iota_\mathcal{H}d(\bar{z})\iota_\mathcal{H})\cdot
(\iota_\mathcal{H}d(\bar{z})^{-1}\iota_\mathcal{H})=I_\mathcal{H},$$
one has $\iota_\mathcal{H}d(\bar{z})^{-1}\iota_\mathcal{H}=
(\iota_\mathcal{H}d(\bar{z})\iota_\mathcal{H})^{-1}=d(z)^{-1}$.
Therefore,
\begin{eqnarray*}
f^\sharp(z)&=&
\iota_\mathcal{U}f(\bar{z})\iota_\mathcal{U}=\iota_\mathcal{U}
(a(\bar{z})-b(\bar{z})d(\bar{z})^{-1} c(\bar{z}))\iota_\mathcal{U}
\\
&=&  \iota_\mathcal{U}a(\bar{z})\iota_\mathcal{U}-(
\iota_\mathcal{U}b(\bar{z})\iota_\mathcal{H})\cdot
(\iota_\mathcal{H}d(\bar{z})^{-1}\iota_\mathcal{H})\cdot
(\iota_\mathcal{H}c(\bar{z})\iota_\mathcal{U}) \\
&=& a(z)-b(z)d(z)^{-1}c(z)=f(z),
\end{eqnarray*}
i.e., $f$ is $\iota_\mathcal{U}$-real. Thus, (i) follows from
(ii).

(iv)$\Rightarrow$(i). Let (iv) hold. Then the operator
$U=U^*=U^{-1}$ is
$(\iota_\mathcal{X}\oplus\iota_\mathcal{U})$-real (and, by the
way, $(\iota_\mathcal{X}\oplus\iota_\mathcal{U})$-symmetric due to
Lemma~\ref{lem:iota}), i.e.,
$$\left[\begin{array}{cc}
  \iota_\mathcal{X} & 0 \\
  0 & \iota_\mathcal{U}
\end{array}\right]\left[\begin{array}{cc}
  A & B \\
  C & D
\end{array}\right]\left[\begin{array}{cc}
  \iota_\mathcal{X} & 0 \\
  0 & \iota_\mathcal{U}
\end{array}\right]=\left[\begin{array}{cc}
  A & B \\
  C & D
\end{array}\right].$$
This is equivalent to the following identities:
$$\iota_\mathcal{X}A\iota_\mathcal{X}=A,\quad \iota_\mathcal{X}B\iota_\mathcal{U}=B,
\quad \iota_\mathcal{U}C\iota_\mathcal{X}=C,\quad
\iota_\mathcal{U}D\iota_\mathcal{U}=D.$$ Moreover, since
$\iota_\mathcal{X}$ commutes with $P_k,\ k=1,\ldots,N$, one has
$\iota_\mathcal{X}(I_\mathcal{X}-AP(w))\iota_\mathcal{X}=I_\mathcal{X}-AP(\bar{w})$,
and $\iota_\mathcal{X}(I_\mathcal{X}-AP(w))^{-1}\iota_\mathcal{X}=
(I_\mathcal{X}-AP(\bar{w}))^{-1}$ (we already used an analogous
argument above). Therefore,
\begin{eqnarray*}
\mathcal{F}^\sharp(w) &=&
\iota_\mathcal{U}\mathcal{F}(\bar{w})\iota_\mathcal{U}=
\iota_\mathcal{U}[D+CP(\bar{w})(I_\mathcal{X}-AP(\bar{w}))^{-1}B]\iota_\mathcal{U}
\\ &=& D+CP(w)(I_\mathcal{X}-AP(w))^{-1}B=\mathcal{F}(w),\quad
w\in\mathbb{D}^N,
\end{eqnarray*}
i.e., $\mathcal{F}$ is $\iota_\mathcal{U}$-real. Applying the
inverse double Cayley transform to $\mathcal{F}$, one can see that
$f$ is also $\iota_\mathcal{U}$-real on $\Pi^N$, and hence, on
$\Omega_N$. Thus, (i) follows from (iv).

(iii)$\Rightarrow$(iv). Let $f$ satisfy the identity
(\ref{eq:kerid}) with holomorphic $\iota_\mathcal{U}$-real PSD
kernels $\Phi_k(z,\zeta),\ k=1,\ldots,N$, on $\Pi^N\times\Pi^N$.
Arguing like in the proof of (iii)$\Rightarrow$(ii) above, we get
Hilbert spaces $\mathcal{M}_k$, holomorphic
$L(\mathcal{U,M}_k)$-valued functions $\varphi_k$ on $\Pi^N$, such
that $\Phi_k(z,\zeta)=\varphi_k(\zeta)^*\varphi_k(z),\
(z,\zeta)\in\Pi^N\times\Pi^N$, AUIs $\iota_{\mathcal{M}_k}$ on
$\mathcal{M}_k,\ k=1,\ldots,N$, and
$\iota_\mathcal{M}=\bigoplus_{k=1}^N\iota_{\mathcal{M}_k}$ on
$\mathcal{M}=\bigoplus_{k=1}^N\mathcal{M}_k$, for which
(\ref{eq:iota-m'}) holds, and hence
$\iota_\mathcal{M}\varphi(z)=\varphi(\bar{z})\iota_\mathcal{U},\
z\in\Pi^N$. Following the proof of Theorem~\ref{thm:double-C}, we
get consecutively: identities (\ref{eq:AH+}) and  (\ref{eq:AH-})
with holomorphic PSD kernels $\Xi_k(w,\omega)$ such that
$\Xi_k(w,\omega)=\xi_k(\omega)^*\xi_k(w),\
(w,\omega)\in\mathbb{D}^N\times\mathbb{D}^N$, and $\xi_k,\
k=1,\ldots,N$, are given in (\ref{eq:xi}), moreover
$\iota_\mathcal{M}\xi(w)=\xi(\bar{w})\iota_\mathcal{U},\
w\in\mathbb{D}^N$; then identities (\ref{eq:Ag+}) and
(\ref{eq:Ag-}) with holomorphic PSD kernels $\Theta_k(w,\omega)$
such that $\Theta_k(w,\omega)=\theta_k(\omega)^*\theta_k(w),\
(w,\omega)\in\mathbb{D}^N\times\mathbb{D}^N$, and $\theta_k,\
k=1,\ldots,N$, are given in (\ref{eq:theta}), moreover
$\iota_\mathcal{M}\theta(w)=\theta(\bar{w})\iota_\mathcal{U},\
w\in\mathbb{D}^N$; then identities (\ref{eq:Ag++}) and
(\ref{eq:Ag--}) with holomorphic functions
$\widetilde{\theta_k}(w)=\theta_k(\bar{w})^*,\ w\in\mathbb{D}^N$,
taking values in $L(\mathcal{M}_k,\mathcal{U})$, moreover,
$\iota_\mathcal{U}\widetilde{\theta_k}(w)=\widetilde{\theta_k}
(\bar{w})\iota_{\mathcal{M}_k},\ w\in\mathbb{D}^N$. The latter
equality is valid since for any $m\in\mathcal{M}_k,\
 u\in\mathcal{U}$ one has
\begin{eqnarray*}
\langle\iota_\mathcal{U}\widetilde{\theta_k}(w)m,u\rangle &=&
\langle\iota_\mathcal{U}\theta_k(\bar{w})^*m,u\rangle
=\langle\iota_\mathcal{U}u,\theta_k(\bar{w})^*m\rangle =
\langle\theta_k(\bar{w})\iota_\mathcal{U}u,m\rangle \\
 =\langle\iota_{\mathcal{M}_k}\theta_k(w)u,m\rangle
&=& \langle\iota_{\mathcal{M}_k}m,\theta_k(w)u\rangle =
\langle\theta_k(w)^*\iota_{\mathcal{M}_k}m,u\rangle =
\langle\widetilde{\theta_k}(\bar{w})\iota_{\mathcal{M}_k}m,u\rangle.
\end{eqnarray*}
Denote $\widetilde{\theta}(w):=[\begin{array}{ccc}
  \widetilde{\theta_1}(w) & \ldots & \widetilde{\theta_N}(w)
\end{array}]$. Then $\iota_\mathcal{U}\widetilde{\theta}(w)=\widetilde{\theta}(\bar{w})
\iota_\mathcal{M},\ w\in\mathbb{D}^N$. For the reproducing kernels
$K_k(w,\omega)=\widetilde{\theta}(w)\widetilde{\theta}(\omega)^*$
of the spaces $\mathcal{M}_k$ we get the identities
$\iota_\mathcal{U}K_k(\bar{w},\bar{\omega})\iota_\mathcal{U}=K_k(w,\omega),\
(w,\omega)\in\mathbb{D}^N\times\mathbb{D}^N,\ k=1,\ldots,N$.
Define in the space
$\mathcal{X}=\bigoplus_{k=1}^N\mathcal{D}(\widehat{K_k})$, where
$\mathcal{D}(\widehat{K_k}),\ k=1,\ldots,N$, are Hilbert spaces
with the reproducing kernels
$$\widehat{K_k}(w,\omega)=\left[\begin{array}{cc}
  K_k(w,\omega) & K_k(w,\omega) \\
  K_k(w,\omega) & K_k(w,\omega)
\end{array}\right],$$
the operator $\iota_\mathcal{X}:=\bigoplus_{k=1}^N
(\iota_{\mathcal{M}_k}\oplus\iota_{\mathcal{M}_k})$. Then
$\iota_\mathcal{X}$ is an AUI on $\mathcal{X}$, and
$\iota_\mathcal{X}P_k=P_k\iota_\mathcal{X},\ k=1,\ldots,N$.
Moreover, it is easy to see that
$(\iota_\mathcal{X}\oplus\iota_\mathcal{U})U_0=
U_0(\iota_\mathcal{X}\oplus\iota_\mathcal{U})$, and therefore
$(\iota_\mathcal{X}\oplus\iota_\mathcal{U})\widetilde{U_0}=
\widetilde{U_0}(\iota_\mathcal{X}\oplus\iota_\mathcal{U})$. It is
clear that $\mathcal{D}_0$ is invariant under
$\iota_\mathcal{X}\oplus\iota_\mathcal{U}$, as well as
$\mathcal{D}_0^\perp=(\mathcal{X}\oplus\mathcal{U})\ominus\overline{\mathcal{D}_0}$
(in fact, since $\mathcal{D}_0\supset\mathcal{U}$, one has
$(\mathcal{X}\oplus\mathcal{U})\ominus\overline{\mathcal{D}_0}\subset\mathcal{X}$,
and $\mathcal{D}_0^\perp$ is invariant under $\iota_\mathcal{X}$).
Indeed, for any $h_1\in\mathcal{D}_0,\ h_2\in\mathcal{D}_0^\perp$
one has
$$\langle(\iota_\mathcal{X}\oplus\iota_\mathcal{U})h_2,h_1\rangle=
\langle(\iota_\mathcal{X}\oplus\iota_\mathcal{U})h_1,h_2\rangle=0,$$
since
$(\iota_\mathcal{X}\oplus\iota_\mathcal{U})h_1\in\mathcal{D}_0$,
thus
$(\iota_\mathcal{X}\oplus\iota_\mathcal{U})h_2\in\mathcal{D}_0^\perp$.
As $U=\widetilde{U}_0\oplus I_{\mathcal{D}_0^\perp}$, we get
$(\iota_\mathcal{X}\oplus\iota_\mathcal{U})U=
U(\iota_\mathcal{X}\oplus\iota_\mathcal{U})$. So, we see that (iv)
follows from (iii). The proof is complete.
\end{proof}

\section{Conclusion and open problems}\label{sec:conclusion}
In this paper the class $\mathbb{R}\mathcal{B}_N^{n\times n}$,
which was defined by M.~F.~Bessmertny\u{\i} in \cite{Bes} (see
also \cite{Bes1}) as a class of rational $n\times n$ matrix-valued
functions having a long resolvent representation (\ref{eq:sc})
with matrix coefficients $A_0=0,\ A_k=A_k^*=A_k^T\geq 0,\
k=1,\ldots,N$, in (\ref{eq:lp}) (note, that matrices $A_k$ have
real entries), was generalized in several directions
simultaneously. First, one can consider the class
$\mathcal{B}_N^{n\times n}=\mathbb{C}\mathcal{B}_N^{n\times n}$ of
rational $n\times n$ matrix-valued functions having a long
resolvent representation (\ref{eq:sc}) with matrix coefficients
$A_0=0,\ A_k=A_k^*\geq 0,\ k=1,\ldots,N$, in (\ref{eq:lp}) (i.e.,
entries of matrices $A_k$ are complex, not necessarily real).
Second, we have introduced the class $\mathcal{B}_N(\mathcal{U})$,
which is a generalization of $\mathcal{B}_N^{n\times n}$,
consisting of holomorphic (not necessarily rational) functions on
the open right polyhalfplane $\Pi^N$ (and naturally extendable to
the domain
$\Omega_N=\bigcup_{\lambda\in\mathbb{T}}(\lambda\Pi)^N$) which
take values in $L(\mathcal{U})$ for a (not necessarily
finite-dimensional) Hilbert space $\mathcal{U}$ and having
representations of the form (\ref{eq:sc})  with coefficients
$A_0=0,\ A_k=A_k^*\geq 0,\ k=1,\ldots,N$, in (\ref{eq:lp}), which
are linear bounded operators on $\mathcal{U}\oplus\mathcal{H}$,
where a Hilbert space $\mathcal{H}$ is not supposed to be
finite-dimensional. We have obtained several equivalent
characterizations of the class $\mathcal{B}_N(\mathcal{U})$, which
we call the Bessmertny\u{\i} class, scattered in different parts
of this paper, and for convenience of a reader we collect them now
in the following theorem.
\begin{thm}\label{thm:fin}
Let $f$ be a holomorphic function on the domain
$\Omega_N=\bigcup_{\lambda\in\mathbb{T}}(\lambda\Pi)^N\subset\mathbb{C}^N$
which takes values in $L(\mathcal{U})$ for a Hilbert space
$\mathcal{U}$. Then the following statements are equivalent:
\begin{description}
    \item[(i)] There exist a Hilbert space
$\mathcal{H}$ and a representation
$$f(z)=a(z)-b(z)d(z)^{-1}c(z),\quad z\in\Omega_N,$$
of $f$, where
$$A(z)=z_1A_1+\cdots +z_NA_N=\left[\begin{array}{cc}
  a(z) & b(z) \\
  c(z) & d(z)
\end{array}\right]\in L(\mathcal{U}\oplus\mathcal{H}),$$
with $A_k=A_k^*\geq 0,\ k=1,\ldots,N$;
    \item[(ii)] there exists a representation
$$f(z)=\sum_{k=1}^Nz_k\Phi_k(z,\zeta),\quad
(z,\zeta)\in\Omega_N\times\Omega_N,$$ of $f$, where
$\Phi_k(z,\zeta),\ k=1,\ldots,N$, are holomorphic PSD kernels on
$\Omega_N\times\Omega_N$;
    \item[(iii)] $f$ satisfies the conditions:
    \begin{enumerate}
        \item $f(\lambda z_1,\ldots,\lambda z_N)=\lambda f(z_1,\ldots,z_N),\
        \lambda\in\mathbb{C}\backslash\{ 0\},\
        z\in\Omega_N$;
        \item $f(\mathbf{R})+f(\mathbf{R})^*\geq 0,\quad
        \mathbf{R}\in\mathcal{A}^N$\\
    (the set $\mathcal{A}^N$ of $N$-tuples of commuting strictly
    accretive operators on a Hilbert space, and the functional
    calculus for $f(\mathbf{R})$ are defined in
    Section~\ref{sec:b-calc});
        \item $f(\bar{z})=f(z)^*,\quad z\in\Omega_N$;
    \end{enumerate}
    \item[(iv)] there exist Hilbert spaces
        $\mathcal{X},\mathcal{X}_1,\ldots,\mathcal{X}_N$, such
        that $\mathcal{X}=\bigoplus_{k=1}^N\mathcal{X}_k$, and an
        Agler representation
        $$\mathcal{F}(w)=D+CP(w)(I_\mathcal{X}-AP(w))^{-1}B,$$
    of a double Cayley transform $\mathcal{F}=\mathcal{C}(f)$ of $f$
    (which is defined by (\ref{eq:Cayley-val})),
     where $P(w)=\sum_{k=1}^Nw_kP_{\mathcal{X}_k}$, and
     $$\left[\begin{array}{cc}
       A & B \\
       C & D
     \end{array}\right]=U=U^{-1}=U^*\in
     L(\mathcal{X}\oplus\mathcal{U}).$$
\end{description}
\end{thm}
 Thus,
 $\mathcal{B}_N(\mathcal{U})$ can be defined as a class of
 functions $f$ satisfying any (and hence, all)
of conditions (i)--(iv) of Theorem~\ref{thm:fin}.

We have introduced also the class
$\iota\mathbb{R}\mathcal{B}_N(\mathcal{U})$ consisting of all
$\iota$-real functions from  $\mathcal{B}_N(\mathcal{U})$, for an
anti-unitary involution $\iota=\iota_\mathcal{U}$ on
$\mathcal{U}$. The operator $\iota$ plays a role analogous to the
complex conjugation on $\mathbb{C}^N$, and the class
$\iota\mathbb{R}\mathcal{B}_N(\mathcal{U})$ is a generalization of
the class $\mathbb{R}\mathcal{B}_N^{n\times n}$. We have obtained
several characterizations of the class
$\iota\mathbb{R}\mathcal{B}_N(\mathcal{U})$, which are collected
in Theorem~\ref{thm:real}.

Let us note that though several descriptions of the classes
$\mathcal{B}_N(\mathcal{U})$ and
$\iota\mathbb{R}\mathcal{B}_N(\mathcal{U})$ were obtained in this
paper, the investigation of these classes is still far from its
final point.  We formulate and discuss below the most important
questions subject to further investigation.

Recall that the class $\mathcal{B}_N(\mathcal{U})$ (resp.,
$\iota\mathbb{R}\mathcal{B}_N(\mathcal{U})$) is a subclass of
$\mathcal{P}_N(\mathcal{U})$ (resp.,
$\iota\mathbb{R}\mathcal{P}_N(\mathcal{U})$), where the latter is
a class of all $L(\mathcal{U})$-valued functions holomorphic in
$\Omega_N$ and satisfying the conditions:
\begin{enumerate}
    \item $f(\lambda z_1,\ldots,\lambda z_N)=\lambda f(z_1,\ldots,z_N),\
        \lambda\in\mathbb{C}\backslash\{ 0\},\
        z\in\Omega_N$;
        \item $f(z)+f(z)^*\geq 0,\quad
        z\in\Pi^N$;
        \item $f(\bar{z})=f(z)^*\quad (=f^\sharp(z):=\iota f(\bar{z})\iota),
        \quad z\in\Omega_N$.
\end{enumerate}
Analogously, the class  $\mathcal{B}_N^{n\times n}$ (resp.,
$\mathbb{R}\mathcal{B}_N^{n\times n}$) is a subclass of
$\mathcal{P}_N^{n\times n}$ (resp.,
$\mathbb{R}\mathcal{P}_N^{n\times n}$), where the latter is a
class of all rational $n\times n$ matrix-valued functions
satisfying the conditions:
\begin{enumerate}
    \item $f(\lambda z_1,\ldots,\lambda z_N)=\lambda f(z_1,\ldots,z_N),\
        \lambda\in\mathbb{C}\backslash\{ 0\},\
        z\in\mathbb{C}^N$;
        \item $f(z)+f(z)^*\geq 0,\quad
        z\in\Pi^N$;
        \item $f(\bar{z})=f(z)^*\quad (=f(\bar{z})^T),
        \quad z\in\mathbb{C}^N$.
\end{enumerate}
It is known that  in the cases $N=1$ and $N=2$ condition~2 for all
of these classes is equivalent to condition (iii2) of
Theorem~\ref{thm:fin}, thus for these cases the class
$\mathcal{B}_N(\mathcal{U})$ (resp.,
$\iota\mathbb{R}\mathcal{B}_N(\mathcal{U})$,
$\mathcal{B}_N^{n\times n}$, $\mathbb{R}\mathcal{B}_N^{n\times
n}$) coincides with $\mathcal{P}_N(\mathcal{U})$ (resp.,
$\iota\mathbb{R}\mathcal{P}_N(\mathcal{U})$,
$\mathcal{P}_N^{n\times n}$, $\mathbb{R}\mathcal{P}_N^{n\times
n}$).
\begin{prob}
For which $N\geq 3,\ \mathcal{U},\ \iota_\mathcal{U}$ and
$n\in\mathbb{N}$ the class $\mathcal{B}_N(\mathcal{U})$ (resp.,
$\iota\mathbb{R}\mathcal{B}_N(\mathcal{U})$,
$\mathcal{B}_N^{n\times n}$, $\mathbb{R}\mathcal{B}_N^{n\times
n}$) is a proper subclass of $\mathcal{P}_N(\mathcal{U})$ (resp.,
$\iota\mathbb{R}\mathcal{P}_N(\mathcal{U})$,
$\mathcal{P}_N^{n\times n}$, $\mathbb{R}\mathcal{P}_N^{n\times
n}$)?
\end{prob}
For the classes $\mathcal{B}_N(\mathcal{U})$ and
$\iota\mathbb{R}\mathcal{B}_N(\mathcal{U}),\ N\geq 3$, this
problem, can be reformulated as follows: find
$f\in\mathcal{P}_N(\mathcal{U})$ and $\mathbf{R}\in\mathcal{A}^N$
such that the selfadjoint operator $f(\mathbf{R})+f(\mathbf{R})^*$
is not PSD. The latter is equivalent to the inequality
$\|\mathcal{F}(\mathbf{T})\| >1$, where
$\mathcal{F}=\mathcal{C}(f)$ is a holomorphic contractive
$L(\mathcal{U})$-valued function on $\mathbb{D}^N$, and the
$N$-tuple of operators
$\mathbf{T}=(T_1,\ldots,T_N)\in\mathcal{C}^N$ is defined by
$T_k:=(R_k-I)(R_k+I)^{-1},\ k=1,\ldots,N$, i.e. for $\mathcal{F}$
and $\mathbf{T}$ the \emph{generalized von Neumann inequality}
fails. There are examples of holomorphic contractive
operator-valued functions $\mathcal{F}$ and $N$-tuples
$\mathbf{T}$ of commuting strict contractions on a Hilbert space,
for which the generalized von Neumann inequality fails (see
\cite{V,CD,K1}), however in these examples the requirement
$\mathcal{F}=\mathcal{C}(f)$ for some
$f\in\mathcal{P}_N(\mathcal{U})$ is not fulfilled. A function
$\mathcal{F}$ satisfying this requirement must have a certain
complicated structure induced by the homogeneity structure of $f$.
Thus, more sophisticated examples should be found to meet this
condition.

Another open problem concerns to characterization of classes
$\mathcal{B}_N^{n\times n}$ and $\mathbb{R}\mathcal{B}_N^{n\times
n}$ (the formulation below is given for the first of them).
\begin{prob}\label{prob:2}
Does the representation $f(z)=\sum_{k=1}^Nz_k\Phi_k(z,\zeta)$ of
an arbitrary $f\in\mathcal{P}_N^{n\times n}$, where
$\Phi_k(z,\zeta),\ k=1,\ldots,N$, are rational $n\times n$
matrix-valued functions which are holomorphic PSD kernels on
$\Omega_N\times\Omega_N, \ N\geq 3$, imply
$f\in\mathcal{B}_N^{n\times n}$?
\end{prob}
Let us formulate this more accurately.
\begin{prob}
For which $N\geq 3,\ n\in\mathbb{N}$, and
$f\in\mathcal{P}_N^{n\times n}$ the question in
Problem~\ref{prob:2} has a positive answer?
\end{prob}




\end{document}